\documentclass[preprint,12pt]{article} 
\usepackage{b}
\usepackage[a4paper]{geometry}
\usepackage[small,it]{caption}
\usepackage{empheq}

\usepackage{subcaption}

\usepackage{cleveref}
\crefname{equation}{}{}
\Crefname{equation}{}{}
\crefname{eqnarray}{}{}
\Crefname{eqnarray}{}{}

\def\x{{\bf x}}
\def\y{{\bf y}}
\def\s{{\bf s}}
\def\t{{\bf t}}
\def\bigo{\mathcal{O}}
\def\littleo{o}

\newtheorem{remark}{Remark}
\newtheorem{proposition}{Proposition}

\title{A stabilized separation of variables method for 
the modified biharmonic equation}
\author{Travis Askham}
\date{\today}

\begin{document}

\maketitle


\begin{abstract}
  The modified biharmonic equation is encountered in a variety
  of application areas, including streamfunction formulations of 
  the Navier-Stokes equations. We develop a separation of variables
  representation for this equation in polar coordinates, for either the
  interior or exterior of a disk, and derive a new class of special functions
  which makes the approach stable. We discuss how these functions can be 
  used in conjunction with fast algorithms to accelerate the solution of
  the modified biharmonic equation or the ``bi-Helmholtz" equation in
  more complex geometries.
\end{abstract}



\section{Introduction} 

Many fourth order elliptic partial 
differential equations of physical interest
can be expressed in terms of the composition of
two second order elliptic differential operators. We focus here on the
modified biharmonic equation in two dimensions, but 
much of what follows can be applied to other equations,
such as the bi-Helmholtz equation, and has a natural extension to the three
dimensionsal setting. In a domain $\Omega$
with boundary $\Gamma$,
the modified biharmonic equation can be written in the form 
\begin{equation}\label{eq:modbh}
  (\Delta^2 - \lambda^2 \Delta) u = 
  \Delta( \Delta - \lambda^2 ) u = 0 \mbox{ in } \Omega \; ,
\end{equation}
with $\lambda \in \R$, subject to two application-specific, boundary conditions.
This equation arises naturally when solving the 
the Navier-Stokes equations using an implicit marching scheme
\cite{chapko1997rothe,greengard1998integral,chapko2001combination,
biros2004embedded,Jiang2013}.

Since the Green's function for the governing equation is 
\cite{Jiang2013}
\[ G(r) = - \frac{1}{2 \pi \lambda^2} \left( \ln r + K_0(\lambda r) \right), \]
any solution to \eqref{eq:modbh} can be expressed using separation of variables 
in the interior of a disk by
\begin{equation}
 u(r,\theta) = \sum_{n=-\infty}^\infty
  \left[ \alpha_n r^{|n|} + \beta_n I_n(\lambda r) \right] \, e^{i n \theta}.
\end{equation}
In the exterior of a disk, assuming the solution is bounded,
the general solution takes the form
\begin{equation}
 u(r,\theta) = 
C_0 + C_1 K_0(\lambda r) + \, 
\sum_{\substack{n=-\infty \\ n \neq 0}}^\infty
  \left[ \alpha_n r^{-|n|} + \beta_n K_n(\lambda r) \right] \, e^{i n \theta}.
\end{equation}
Unfortunately, using this representation naively for interior or exterior
boundary value problems on a disk of radius $R$ leads to numerical instabilities
when the value $\lambda R$ is small. While this problem is of
interest in its own right, such expansions play a role in more general domains
when mulitpole  and local expansions are used to represent outgoing and incoming
fields in more general geometries, as explained below.
Here, we propose a stabilized separation of variables
approach based on a new class of special functions. 

\begin{remark} \label{rmk:bihelm}
In the case of the bi-Helmholtz equation, 
\begin{equation}\label{eq:bihelm}
  (\Delta - \lambda_1^2)(\Delta - \lambda_2^2) u = 0,
\end{equation}
the separation of variables representation in the interior of a disk is given by
\begin{equation}
 u(r,\theta) = \sum_{n=-\infty}^\infty
  \left[ \alpha_n I_n(\lambda_1 r) + \beta_n I_n(\lambda_2 r) \right] \, e^{i n \theta}.
\end{equation}
In this setting, it is perhaps clearer that there are two dimensionless
quantities involved: $\lambda_1 R$ and 
$\lambda_2 R$, with obvious ill-conditioning involved when  
$\lambda_1 \approx \lambda_2$.
\end{remark}

The remainder of this paper is
organized as follows. In sections \ref{sec:mathprelim},
and \ref{sec:analysis},
we provide some mathematical preliminaries,
review the classical separation of variables approach
to the modified biharmonic problem, and
define new functions $Q_n$ and $P_n$ for stably
representing solutions in the small $\lambda R$ regime.
In \cref{sec:catastrophe},
we present numerical examples to illustrate the necessity for 
stabilization and to demonstrate the efficacy of
our new functions. In the discussion of
\cref{sec:discuss}, we outline how these functions
can be used for stably solving \cref{eq:modbh}
on more complex geometries with an accelerated integral 
equation method. Such methods will depend on translation
operators for series using the $Q_n$ and $P_n$ functions, 
which we include in \cref{sec:props}.
\Cref{sec:prelim} provides some needed
properties of Bessel and Laurent series. 

\section{Mathematical preliminaries}
\label{sec:mathprelim}

\subsection{Notation}

In the following, we use the ``big O''
notation to describe the order of the remaining
terms in a power series. The expression
$f(\epsilon) = \bigo (g(\epsilon))$ implies that
there exist positive constants $C$ and $\epsilon_0$
such that

\begin{equation}
  |f(\epsilon)| \leq C g(\epsilon) \; ,
\end{equation}
when $0 < \epsilon < \epsilon_0$. As we are
concerned with the asymptotic behavior with
respect to two variables, we often write
$f(\epsilon,\delta) = \bigo ( g(\epsilon,\delta))$
which implies that there exist positive constants
$C$, $\epsilon_0$, and $\delta_0$ such that 

\begin{equation}
  |f(\epsilon,\delta)| \leq C g(\epsilon,\delta) \; ,
\end{equation}
when $0 < \epsilon < \epsilon_0$ and
$0 < \delta < \delta_0$.

In \cref{sec:discuss}, we use the ``big O''
notation to describe computational cost. For
this case, the expression $f(N) = \bigo(g(N))$
implies that there exist positive constants $C$
and $N_0$ such that

\begin{equation}
  |f(N)| \leq C g(N) \; ,
\end{equation}
when $N > N_0$.

\subsection{Condition numbers of $2\times 2$ linear
  systems and diagonal scaling}

\label{sec:cond2by2}

In this section, we review some basic results
from linear algebra.
Consider an invertible linear system of the form

\begin{equation}
  Ax = b \; .
\end{equation}
The condition number $\kappa(A)$ of the matrix $A$
describes the sensitivity of the problem of
recovering $x$ from $b$ \cite{trefethen1997numerical}.
Suppose that an approximate
solution $x_0$ is found such that

\begin{equation}
  \dfrac{\|b-Ax_0\|_2}{\|b\|_2} = \epsilon \; ,
\end{equation}
where $\|\cdot \|_2$ denotes the Euclidean norm.
Let $\kappa(A) = \sigma_{\mbox{max}}(A)/\sigma_{\mbox{min}}(A)$,
where $\sigma_{\mbox{max}}(A)$ and $\sigma_{\mbox{min}}(A)$ denote
the maximum and minimum singular values of $A$.
Then, $x_0$ satisfies

\begin{equation}
  \dfrac{\|x-x_0\|_2}{\|x\|} \leq \kappa(A) \epsilon \; .
\end{equation}

Suppose that each column $a_i$ of $A$ represents some
function from a basis and that $x$ represents the
coefficients which reconstruct $b$ in that basis.
In this case, the notion of the sensitivity of $x$
to changes in $b$ should be unaffected by scaling
the columns of $A$.
Let $D$ be an invertible diagonal matrix and $A$, $x$,
$b$, and $x_0$ be as above. We note that the residual
is unaffected by scaling $x_0$ and $A$, i.e. that
$\|b-Ax_0\|_2 = \|b-(AD)(D^{-1}x_0)\|_2$ but that the
condition numbers of $A$ and $AD$ can be significantly
different. To avoid this ambiguity, we quantify the
sensitivity of recovering $x$ from $b$ in terms of
the condition number of $\tilde{A} = AD$, where
$D$ is a diagonal matrix with $D_{ii} =1/ \|a_i\|_2$.

For $2 \times 2$ matrices, this is a natural
normalization. It is straightforward to prove

\begin{proposition}
  Let $A$ be a $2\times 2$ matrix with columns $a_i$.
  If $D$ is a diagonal matrix with
  $D_{ii} = 1/\|a_i\|_2$, then

  \begin{equation}
    \kappa (AD) = \min_{v \in \R^2} \kappa (A \diag (v) ) \; ,
  \end{equation}
  where $\diag(v)$ denotes the diagonal matrix
  whose main diagonal is given by $v$.
\end{proposition}

It is particularly simple to characterize the
condition number of matrices with this
scaling. We have

\begin{proposition}
  Let $A$ be a $2\times 2$ matrix with columns
  denoted by $a_i$. Suppose that $\|a_i\| = 1$.
  Then

  \begin{equation}
    \kappa = \sqrt{\dfrac{1+c}{1-c}} \; ,
  \end{equation}
  where $c = |a_1^\intercal a_2|$ is the cosine
  of the angle between the two columns of $A$.
\end{proposition}

\subsection{Separation of variables}
\label{sec:sep}

Consider the modified biharmonic equation \cref{eq:modbh}
where $\Omega$ is
the disk of radius $R$ centered at the origin.
For given functions $f$ and $g$, we prescribe 
Dirichlet boundary conditions on $u$, i.e.

\begin{align}
  u &= f \mbox{ on } \Gamma \; , \\
  \partial_n u &= g \mbox{ on } \Gamma \; ,
\end{align}
where $\partial_n$ denotes the outward normal
derivative. Other types of boundary conditions
may be considered.

To take advantage of the simplicity of this 
geometry, we translate the problem to a polar
coordinate system. Let $(r,\theta)$ denote the 
usual polar coordinates
for the point $(x,y) \in \R^2$ as in the following
change of variables

\begin{equation}
\left \{
\begin{array}{rcl}
  x &=& r\cos(\theta) \\ 
  y &=& r\sin(\theta)
\end{array} \right.
\leftrightarrow \;
\left\{
\begin{array}{rcl}
 r &=& \sqrt{x^2+y^2} \\
\theta &=& \arctan(y/x)
\end{array} \right. \; ,
\end{equation}
with $\theta \in [-\pi,\pi)$. Then, the Laplacian
is given by

\begin{equation}
  \Delta = \partial_{rr} + \dfrac1r \partial_r + \dfrac1{r^2} 
  \partial_{\theta\theta} \; .
\end{equation}

Consider the boundary data as functions of the
coordinate $\theta$.
For $f$ and $g$ sufficiently smooth, the Fourier
series 

\begin{align}
  f(R,\theta) &= \sum_{n=-\infty}^\infty f_n e^{in\theta} \; ,\\
  g(R,\theta) &= \sum_{n=-\infty}^\infty g_n e^{in\theta} \; ,
\end{align}
converge uniformly in $\theta$ and the coefficients $f_n$ and
$g_n$ are given by

\begin{align}
  f_n &= \dfrac1{2\pi} \int_{-\pi}^\pi f(R,\theta) 
  e^{-in\theta} \, d\theta \; , \label{eq:fn} \\
  g_n &= \dfrac1{2\pi} \int_{-\pi}^\pi g(R,\theta) 
  e^{-in\theta} \, d\theta \; \label{eq:gn}.
\end{align}

Following standard practice \cite{Fourier1808,Seeley1966},
we seek a solution $u$ of the form 

\begin{equation}
  u(r,\theta) = \sum_{n=-\infty}^\infty u_n(r) e^{in\theta} \; .
\end{equation}
Plugging this form for $u$ into \cref{eq:modbh}, we obtain,
after some simplification,

\begin{equation}
  \sum_{n=-\infty}^\infty \left ( \partial_{rr} + \dfrac1r \partial_r 
  - \dfrac{n^2}{r^2} \right ) \left ( \partial_{rr} + \dfrac1r 
  \partial_r - \dfrac{n^2}{r^2} - \lambda^2 \right )
  u_n(r) e^{in\theta} = 0 \; .
\end{equation}

The above implies that the radial function $u_n(r)$ satisfies
the following ordinary differential equation

\begin{equation}
  \label{eq:ode}
  \left ( \dfrac{d^2}{dr^2} + \dfrac1r \dfrac{d}{dr} - \dfrac{n^2}{r^2} 
  \right ) \left ( \dfrac{d^2}{dr^2} + \dfrac1r \dfrac{d}{dr} 
  - \dfrac{n^2}{r^2}  - \lambda^2 \right ) u_n(r) = 0 \, ,
\end{equation}
subject to certain boundary or regularity conditions. For $n\neq -1,0,1$, we
have

\begin{equation}
  \begin{array}{rclrcl}
    u_n(R) &=& f_n \; , \qquad & u_n'(R) &=& g_n \; , \\
    u_n(0) &=& 0 \; , \qquad & u_n'(0) &=& 0 \; .
  \end{array}
\end{equation}
For $n = 0$, we have

\begin{equation}
  \begin{array}{rclrcl}
    u_0(R) &=& f_0 \; , \qquad & u_0'(R) &=& g_0 \; , \\
    u_0'(0) &=& 0 \; , \qquad & u_0'''(0) &=& 0 \; .
  \end{array}
\end{equation}
Finally, for $n = -1,1$, we have

\begin{equation}
  \begin{array}{rclrcl}
    u_n(R) &=& f_n \; , \qquad & u_n'(R) &=& g_n \; , \\
    u_n(0) &=& 0 \; , \qquad & u_n''(0) &=& 0 \; .
  \end{array}
\end{equation}
The conditions at $r=0$ are derived by assuming that
$u$ has four continuous derivatives at the origin.

It is well known \cite{Abramowitz1966,Olver2010}
that equation \cref{eq:ode} has the following four 
linearly independent solutions for $n\neq 0$: 
$r^{|n|},r^{-|n|},I_n(\lambda r), K_n(\lambda r)$,
where $I_n$ and $K_n$ are the modified Bessel functions.
For $n=0$, the functions $1, \log r, I_0(\lambda r), 
K_0(\lambda r)$ are linearly independent solutions.
The regularity of the solution at zero eliminates
$K_n(\lambda r), r^{-|n|}, \log r$ from the acceptable 
solution set. Therefore, the allowed functions $u_n$ are 
linear combinations of the following form:

\begin{equation}
  u_n(r) = \alpha_n r^{|n|} + \beta_n I_n(\lambda r) \; .
\end{equation}

The boundary conditions for $u_n$ determine 
$\alpha_n$ and $\beta_n$, with the conditions at 
$r=0$ automatically satisfied. From the 
conditions at $r=R$, we obtain the following
linear system for $\alpha_n$ and $\beta_n$

\begin{equation}
\label{eq:coeffs}
\begin{pmatrix}
R^{|n|} & I_n(\lambda R) \\
|n|R^{|n|-1} & 
\dfrac{\lambda}{2} \left (I_{n-1}(\lambda R) + I_{n+1}(\lambda R) \right)
\end{pmatrix}
\begin{pmatrix}
\alpha_n \\ \beta_n
\end{pmatrix}
=
\begin{pmatrix}
f_n \\ g_n
\end{pmatrix} \; .
\end{equation}
The determinant of the system in \cref{eq:coeffs}
is $\lambda R^{|n|} I_{n+1}(\lambda R)$, which is nonzero
for positive $R$. Therefore, the coefficients
$\alpha_n$ and $\beta_n$ are determined by 
the boundary conditions and the above provides
an algorithm for computing $u$. 

In the exterior of a disk, the derivation is analogous
to the above. For simplicity, we consider solutions of 
\cref{eq:modbh} which are bounded with derivatives
that are $\littleo (1/r)$ as $r$ goes to infinity.
This is sufficient for the solutions of the exterior
problem to be unique; see, for example, Proposition 3.5
in \cite{Jiang2013}.
The functions appropriate for the 
exterior of a disk are then $K_n(\lambda r)$ and $r^{-|n|}$.
As in \cref{eq:coeffs}, we obtain a linear system 
for the expansion coefficients

\begin{equation}
\label{eq:coeffs2}
\begin{pmatrix}
R^{-|n|} & K_n(\lambda R) \\
-|n|R^{-|n|-1} & 
-\dfrac{\lambda}{2} \left (K_{n-1}(\lambda R) + K_{n+1}(\lambda R) \right)
\end{pmatrix}
\begin{pmatrix}
\alpha_n \\ \beta_n
\end{pmatrix}
=
\begin{pmatrix}
f_n \\ g_n
\end{pmatrix} \; .
\end{equation}
The determinant of this system is
$-\lambda K_{n-1}(\lambda R)/R^{|n|}$ so that
it is invertible.

In order to make a numerical method out of the
above, one simply truncates the Fourier
series expansions at some finite $N$, i.e.

\begin{align} 
  f(R,\theta) &\approx \sum_{n=-N}^{N+1} f_n e^{in\theta} 
  \; , \\
  g(R,\theta) &\approx \sum_{n=-N}^{N+1} g_n e^{in\theta} 
  \; , \\
  u(r,\theta) &\approx \sum_{n=-N}^{N+1} u_n(r)e^{in\theta}
  \; .
\end{align}
The formulas \cref{eq:fn} and \cref{eq:gn} for 
the coefficients $f_n$ and $g_n$ can be approximated
using the trapezoidal rule with
$M=2N+2$ equispaced points on $\Gamma$ and computed 
rapidly via the fast Fourier transform \cite{Cooley1965}.
The rate of convergence (in $N$) depends on the 
smoothness of the boundary data $f$ and $g$, 
with spectral convergence for analytic $f$ and $g$.

\section{Analysis of the separation of variables
  problem and new basis functions}
\label{sec:analysis}

The difficulty with the above procedure is in 
solving the linear systems \cref{eq:coeffs,eq:coeffs2}.
In particular, for small $\lambda R$, the columns of 
these system matrices are nearly linearly dependent, i.e.,
as $\lambda$ goes to zero, the angle between the
columns goes to zero (there is a similar 
effect for small $R$). As noted in \cref{sec:cond2by2}
this makes the problem of recovering the coefficients, 
$\alpha_n$ and $\beta_n$, from the data, $f_n$
and $g_n$, unstable. In this section, we will 
investigate the nature of this ill-conditioning
and derive new bases which are better conditioned.

We first fix some notation. For a pair of functions 
$(F(r),G(r))$, define the matrix $A(F,G,R)$ to be 

\begin{equation}
  \label{eq:matgeneric}
A(F,G,R) =  \begin{pmatrix}
    F(R) & G(R) \\
    F'(R) & G'(R)
  \end{pmatrix} \; .
\end{equation}
This is the form of the matrix that appears in the 
linear systems \cref{eq:coeffs,eq:coeffs2} used to
solve for the coefficients $\alpha_n$ and $\beta_n$. 
Let $\tilde{B}$ denote the matrix $B$ with its 
columns normalized to unit length. 

For the interior problem, the ill-conditioning of the
basis $(r^{|n|},I_n)$ results from the fact that the $I_n(\lambda r)$ 
and $r^{|n|}$ are very similar functions for small
$r$; they have the same asymptotic behavior to leading order. 
The power series for $I_n(\lambda r)$ is given by 

\begin{equation}
  I_n(\lambda r) = \sum_{k=0}^\infty \dfrac{ \left ( \dfrac{\lambda r}{2} \right )^{2k+|n|} }
  {k! (k+|n|)!} = \dfrac1{2^{|n|} |n|!} \left (\lambda r \right )^{|n|} + \dfrac1{2^{|n|+2} (|n|+1)!}
  \left ( \lambda r \right )^{|n|+2} +  \cdots \; ,
\end{equation}
see \cite[Ch.~10]{NIST:DLMF} for reference. By substituting
this expression into $A(r^{|n|},I_n,R)$, we obtain, for $n \neq 0$,

\begin{equation}
A(r^{|n|},I_n,R) = \begin{pmatrix}
R^{|n|} & R^{|n|} a_n(\lambda) (1+\bigo (\lambda^2R^2)) \\
|n|R^{|n|-1} & |n| R^{|n|-1} a_n(\lambda) (1+\bigo (\lambda^2R^2))
\end{pmatrix} \; ,
\end{equation}
where $a_n(\lambda) = \lambda^{|n|}/(2^{|n|}|n|!)$. We see that
the columns of $\tilde{A}(r^{|n|},I_n,R)$ are nearly co-linear
in the limit as either $\lambda$ or $R$ goes to zero.
The dependence on $\lambda$ and $R$ is not identical (this is wrapped
up in the ``big O'' expressions); in the 
next section, we
see that the condition number of the normalized matrix generally
increases faster as $R$ goes to zero than it does as $\lambda$
goes to zero. For $n=0$, we have the system

\begin{equation}
A(1,I_0,R) =
\begin{pmatrix}
1 & 1+\bigo (\lambda^2R^2) \\
0 & \dfrac12 \lambda^2 R (1+\bigo(\lambda^2 R^2))
\end{pmatrix} \; ,
\end{equation}
so that the condition number increases faster as $\lambda$
goes to zero in this case.

To alleviate this ill-conditioning, we construct basis
functions which are less ``similar" in the
small $R$ and small $\lambda$ limits. 
Let us define the functions $P_n(r)$ by 

\begin{equation} \label{eq:Pn}
  P_n(r) = I_n(\lambda r) 
- \left ( \dfrac{\lambda r}{2} \right )^{|n|} \dfrac1{|n|!} \; ,
\end{equation}
deleting
the first term in the power series for $I_n$. 
Note that, $P_n$ is a solution of \cref{eq:ode} because it is a linear
combination of $I_n$ and $r^{|n|}$.
The matrix we obtain for the basis $(r^{|n|},P_n)$ is

\begin{equation}
A(r^{|n|},P_n,R) = 
\begin{pmatrix}
R^{|n|} & R^{|n|+2} b_n(\lambda) (1+\bigo(\lambda^2 R^2)) \\
|n|R^{|n|-1} & (|n|+2) R^{|n|+1} b_n(\lambda) (1+\bigo(\lambda^2 R^2))
\end{pmatrix} \; ,
\end{equation}
where $b_n(\lambda) = \lambda^{|n|+2}/(2^{|n|+2} (|n|+1)!)$.
As $\lambda$ goes to zero, the columns of $\tilde{A}(r^{|n|},P_n,R)$
do not converge to the same vector, as in the
above. Further,
as $R$ goes to zero, the normalized columns do converge 
to the same limit but
the effect is not as dramatic as for the pair $(r^{|n|},I_n)$.
The problem of recovering the coefficients is
therefore more stable for the basis $(r^{|n|},P_n)$, which
we verify numerically in the next section.

\begin{remark}
  It is simple to evaluate $P_n(r)$ stably. For small 
$\lambda r$,
  the power series for $I_n$, with the first term omitted,
  may be used. For larger $r$, there is no fear of numerical
  cancellation and the formula \cref{eq:Pn} may be used 
  directly, along with existing software for evaluating
  $I_n$. We choose the pair $(r^{|n|},P_n)$ as opposed
  to $(I_n,P_n)$ because $P_n$ and $r^{|n|}$ have different
  asymptotic behavior for large $r$, whereas $P_n$ and
  $I_n$ both grow exponentially.
\end{remark}

\begin{remark}
  If the functions $I_n$ and $r^{|n|}$ were used as a
  basis themselves, their asymptotic similarity in the small $R$ 
  regime would cause other numerical problems for a 
  solution which behaves like $P_n$. 
  To see this, note that for small $r$, the 
  function $P_n(r)$ is $\bigo (r^{|n|+2})$, while the functions 
  $I_n(r)$ and $r^{|n|}$ are $\bigo(r^{|n|})$. Therefore, there 
  is significant numerical cancellation when evaluating $P_n(r)$ 
  via the formula \cref{eq:Pn}. We illustrate this effect
  in the next section.
\end{remark}

A similar analysis applies to the exterior problem. 
The power series for $K_n(\lambda r)$ is given by

\begin{align} \label{eq:Knpow}
  K_{n}\left(\lambda r\right) &=
  \tfrac{1}{2}(\tfrac{1}{2}\lambda r)^{-|n|}
  \sum_{k=0}^{|n|-1}\frac{(|n|-k-1)!}{k!}(-\tfrac{1}{4}\lambda r^{2})^{k}
  +(-1)^{|n|+1}\ln\left(\tfrac{1}{2}\lambda r\right)
  I_{n}\left(\lambda r\right) \nonumber \\
  & \quad +(-1)^{|n|}\tfrac{1}{2}(\tfrac{1}{2}\lambda r)^{|n|}
  \sum_{k=0}^{\infty}\left (\psi\left(k+1\right)+\psi\left(|n|+k+1\right)\right)
  \frac{(\tfrac{1}{4}\lambda r^{2})^{k}}{k!(|n|+k)!} \; .
\end{align}
Substituting this expression into $A(r^{-|n|},K_n,R)$, we
obtain, for $n \neq 0$,

\begin{equation}
A(r^{-|n|},K_n,R) = \begin{pmatrix}
R^{-|n|} & R^{-|n|} c_n(\lambda) (1+\bigo (\lambda^2R^2)) \\
-|n|R^{-|n|-1} & -|n| R^{-|n|-1} c_n(\lambda) (1+\bigo (\lambda^2R^2))
\end{pmatrix} \; ,
\end{equation}
where $c_n(\lambda) = (|n|-1)! 2^{|n|-1} \lambda^{-|n|}$.
Again, after normalization, this linear system is ill-conditioned
as either $R$ or $\lambda$ goes to zero because the normalized columns
become nearly colinear. For the case $n=0$, the basis $(1,K_0)$
results in the system

\begin{equation}
A(1,K_0,R) = \begin{pmatrix}
1 & -\gamma + \log(2) - \log (\lambda R) (1+\bigo (\lambda^2R^2)) \\
0 & -\dfrac1R (1+\bigo (\lambda^2R^2 |\log(\lambda R)|))
\end{pmatrix} \; ,
\end{equation}
which is actually well conditioned, after normalization,
for small $\lambda$ and $R$. For certain exterior problems,
the basis $(\log r, K_0)$
is more appropriate. In this case, we have

\begin{equation}
A(\log r, K_0,R) = \begin{pmatrix}
\log R & -\gamma + \log(2) - \log (\lambda R) (1+\bigo (\lambda^2R^2)) \\
\dfrac1R & -\dfrac1R (1+\bigo (\lambda^2R^2 |\log(\lambda R)|))
\end{pmatrix} \; .
\end{equation}
After normalization, this system is generally well-conditioned
as $\lambda$ goes to zero. As $R$ goes to zero, however,
the two columns become nearly colinear and the normalized
matrix is ill-conditioned. Regardless of these special
cases, the instability for the $n\neq 0$ coefficients
will negatively affect the separation of variables approach.

To avoid this instability, we can define new functions
$Q_n$ for $n\neq 0$ as

\begin{equation}
  \label{eq:Qn}
  Q_n(r) = K_n(\lambda r) - \dfrac{2^{|n|-1} \left (|n|-1\right )!}{\lambda^{|n|}r^{|n|}} \; .
\end{equation}
The function $Q_n$ has a different leading order term from
$K_n$ as $\lambda$ and $r$ go to zero but is still a solution
of \cref{eq:ode} as it is a linear combination of $K_n$ and $r^{-|n|}$.
As noted above, the system \cref{eq:coeffs2} is well conditioned for 
the functions $1$ and $K_0(\lambda r)$. Therefore, the na\"{i}ve
approach works for the zero mode, if you are indeed interested
in solving the exterior problem subject to the conditions we 
have set at infinity. It is convenient, however, to define $Q_0$ as
\begin{equation}
Q_0(r) = K_0(\lambda r) + \log(r) \; .
\end{equation}
This function is closely related to the Green's function for 
the modified biharmonic equation.

As before, the ill-conditioning of the coefficient
recovery problem \cref{eq:coeffs} is improved if
you use the pair of functions $(Q_n,K_n)$ as your
basis. For $|n| > 2$, we obtain a system of the form

\begin{align}
&A(Q_n,K_n,R) = \nonumber \\
& \quad \begin{pmatrix}
  R^{-|n|+2} d_n(\alpha)(1+\bigo (\lambda^2 R^2))
  & R^{-|n|} c_n(\lambda) (1+\bigo (\lambda^2R^2)) \\
  (-|n|+2)R^{-|n|+1} d_n(\alpha) (1+\bigo (\lambda^2 R^2))
 & -|n| R^{-|n|-1} c_n(\lambda) (1+\bigo (\lambda^2R^2))
\end{pmatrix} \; ,
\end{align}
where $d_n(\lambda) = -2^{|n|-3} (|n|-2)! \lambda^{-|n|+2}$.
As $\lambda$ goes to zero, the normalized columns of this
system matrix do not converge to the same vector, as in the
above. Further,
as $R$ goes to zero, the columns do converge but
the effect is not as dramatic as for the pair $(r^{-|n|},K_n)$.
The problem of recovering the coefficients is
therefore more stable for the basis $(Q_n,K_n)$ than it
is for the basis $(r^{-|n|},K_n)$, which
we verify numerically in the next section.

\begin{remark}
  It is simple to evaluate $Q_n(r)$ stably. For small $r$,
  the power series for $K_n$, with the first term omitted,
  may be used. For larger $r$, there is no fear of numerical
  cancellation and the formula \cref{eq:Qn} may be used 
  directly, along with existing software for evaluating
  $K_n$. We choose the pair $(Q_n,K_n)$ as opposed
  to $(r^{-|n|},Q_n)$ because $Q_n$ and $K_n$ have different
  asymptotic behavior for large $r$. 
\end{remark}

\begin{remark}
  If the functions $K_n$ and $r^{-|n|}$ are used as a
  basis, their asymptotic similarity in the small $R$ 
  regime will cause other numerical problems for a 
  solution with terms like $Q_n$. Suppose that $K_n$ and 
  $r^{-|n|}$ are used to evaluate $Q_n$. For small $r$, the 
  function $Q_n(r)$ is $\bigo (r^{-|n|+2})$, while the functions 
  $K_n(r)$ and $r^{-|n|}$ are $\bigo(r^{-|n|})$. Therefore, there 
  is significant numerical cancellation when evaluating $Q_n(r)$ 
  via the formula \cref{eq:Qn}. We demonstrate this effect
  as well in the next section.
\end{remark}

Before proceeding to the numerical experiments, we briefly 
describe the edge cases, i.e. the matrices for $|n| \leq 2$.
When $|n| = 2$, the pair $(Q_n,K_n)$ results in a linear
system of the form

\begin{equation}
A(Q_2,K_2,R)
  \begin{pmatrix}
    -\dfrac12 + \bigo(R^2\lambda^2 |\log(\lambda R)|) &
    \dfrac2{\lambda^2 R^2} + \bigo(1)\\
    -\dfrac14 \lambda^2 R\log(\lambda R) + \bigo(\lambda^2 R)
    & -\dfrac4{\lambda^2R^3} + \bigo(\lambda^2 R |\log(\lambda R)|)
\end{pmatrix} \; .
\end{equation}
The columns of this matrix, after normalization,
do not become colinear as either $\lambda$ or $R$
tends to zero. For the case
that $|n|=1$, the pair $(Q_n,K_n)$ results in a linear
system of the form

\begin{equation}
A(Q_1,K_1,R) = \begin{pmatrix}
  \dfrac12 \lambda R \log(\lambda R) + \bigo (\alpha R) &
  \dfrac{1}{\lambda R} + \bigo( \lambda R |\log(\lambda R)|) \\
  \dfrac12 \lambda \log(\lambda R) + \big(\lambda +R^2) &
  -\dfrac1{\lambda R^2} + \bigo(\lambda |\log(\lambda R)|)
\end{pmatrix} \; ,
\end{equation}
which has nearly colinear columns after normalization
as $R$ goes to zero. The normalized columns are not
colinear in the limit as $\lambda$ goes to zero. Finally,
for the case $n=0$, the pair $(Q_n,K_n)$ results in
a linear system of the form
\begin{equation}
  \begin{pmatrix}
    -\log(\lambda/2) -\gamma + \bigo(\lambda^2R^2 |\log(\lambda R)|
    &
    -\log(\lambda R/2) -\gamma + \bigo(\lambda^2R^2 |\log(\lambda R)| \\
    -\dfrac12 \lambda^2 R \log(\lambda R) + \bigo(\lambda^2 R)
    & -\dfrac1R + \bigo(\lambda^2 R |\log(\lambda R)|)
\end{pmatrix}
\begin{pmatrix}
\alpha_n \\ \beta_n
\end{pmatrix}
=
\begin{pmatrix}
f_n \\ g_n
\end{pmatrix} \; .
\end{equation}
As $\lambda$ goes to zero, the columns of this matrix
slowly become colinear after normalization. As $R$
goes to zero, the normalized columns do not become
colinear.

\section{Numerical tests}
\label{sec:catastrophe}

In this section we present some numerical experiments which
reinforce the ideas of the previous section. The source
code used for these calculations is available online
\cite{askhampapercode}. The Bessel functions were evaluated
using routines from FMMLIB2D \cite{gimbutasfmmlib2d}.
Discrete Fourier transforms were computed using the
FFTPACK \cite{swarztrauber1985fftpack}. To compute
condition numbers, we used the singular value decomposition
routine from EISPACK \cite{garbow1993eispack}.
All code was written in Fortran using double
precision arithmetic and compiled with the gfortran
compiler on Linux.

\subsection{Condition numbers}

For the first test, we verify the analytical observations
about the condition numbers of the matrices described in
the previous section. For a pair of functions $(F(r),G(r))$,
define the matrix $A(F,G,R)$ as in \cref{eq:matgeneric}. 
Let $\tilde{B}$ denote the matrix $B$ with its columns
normalized to unit length. We compute the condition
number of $\tilde{A}(F,G,R)$ where $(F,G)$ is taken
to be each of the bases $(r^{|n|},I_n(\lambda r))$,
$(r^{|n|},P_n(\lambda r))$, $(r^{-|n|},K_n(\lambda r))$, and
$(Q_n(\lambda r),K_n(\lambda r))$ for a range of values
of $n$, $\lambda$, and $R$. To observe the effect of
changing the radius of the domain, we run experiments
with $\lambda$ fixed (at $\lambda = 0.5$) and, for
each $j$ from $-24$ to $8$, ten values of $R$ drawn
uniformly at random from the interval $[2^{j},2^{j+1}]$.
Likewise, to observe the effect of changing the
parameter $\lambda$, we run experiments with $R$ fixed
(at $R = 0.5$) and, for each $j$ from $-24$ to $8$,
ten values of $\lambda$ drawn uniformly at random from
the interval $[2^{j},2^{j+1}]$. In order to compare these,
we plot the condition number as a function of the
product $\lambda R$. 

\begin{figure}[h!]
  \begin{subfigure}[t]{0.5\textwidth}
    \includegraphics[width=\linewidth]{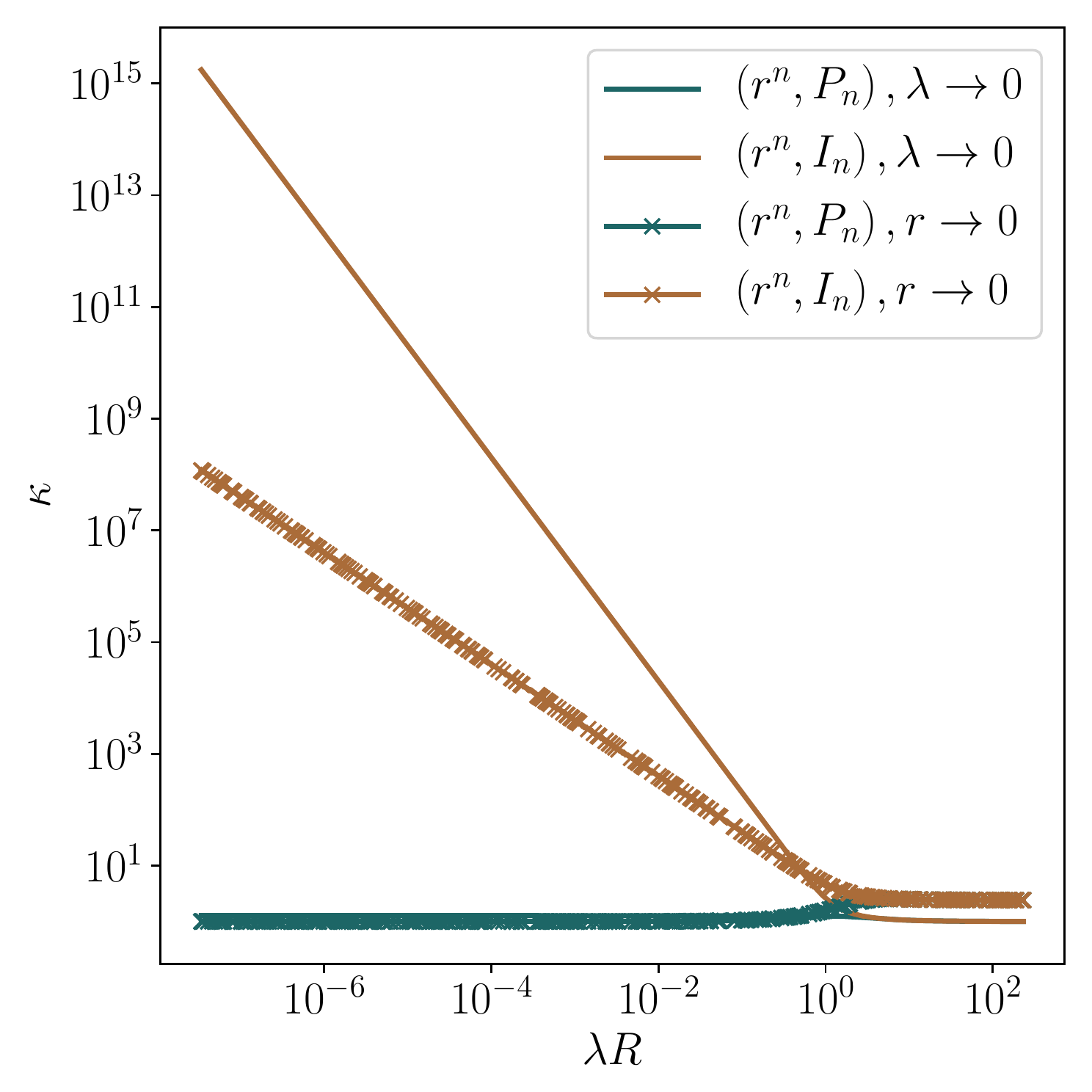}
    \caption{$n = 0$}
  \end{subfigure}
  \begin{subfigure}[t]{0.5\textwidth}
    \includegraphics[width=\linewidth]{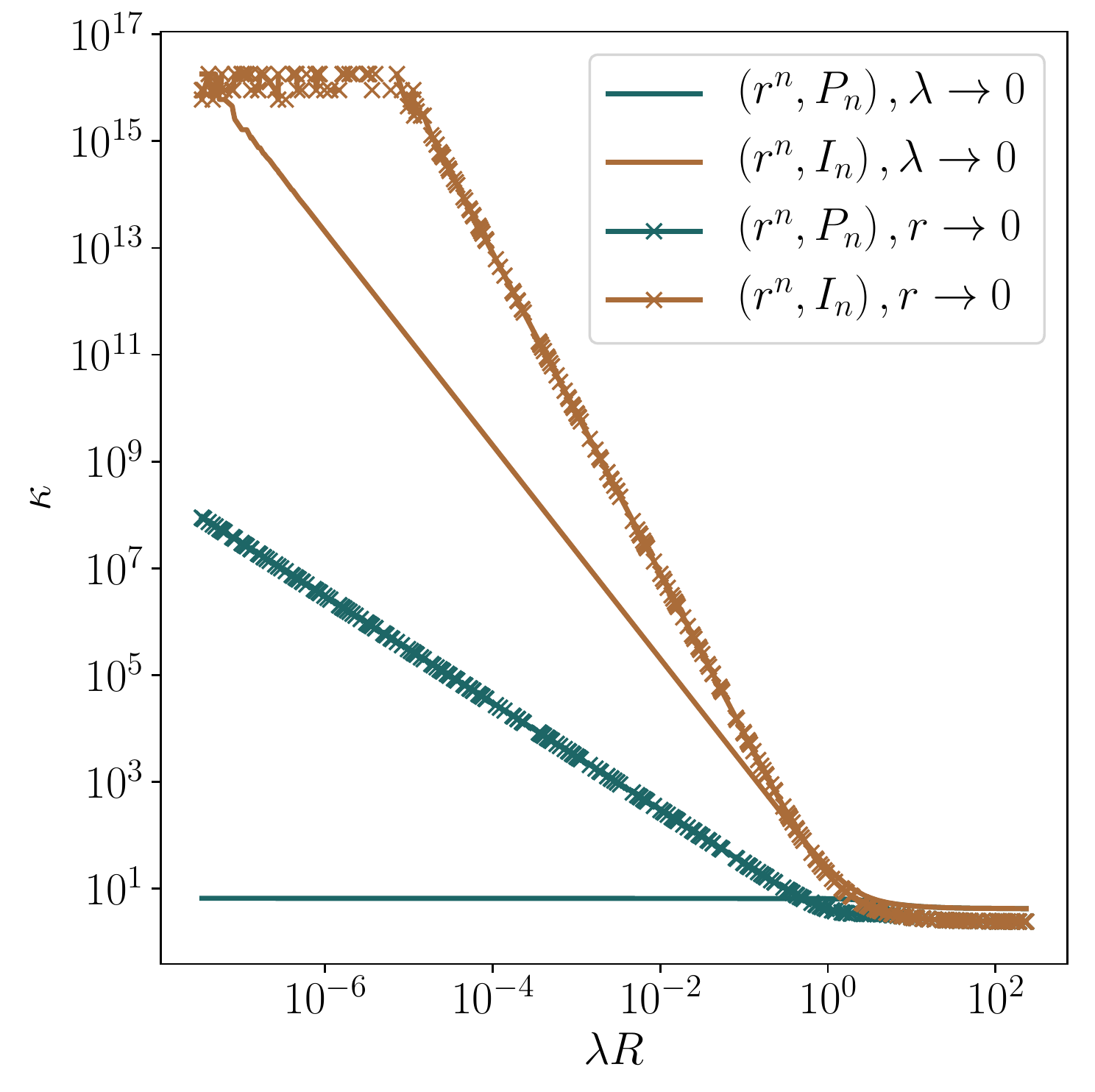}
    \caption{$n=1$}    
  \end{subfigure}
  \begin{subfigure}[t]{0.5\textwidth}
    \includegraphics[width=\linewidth]{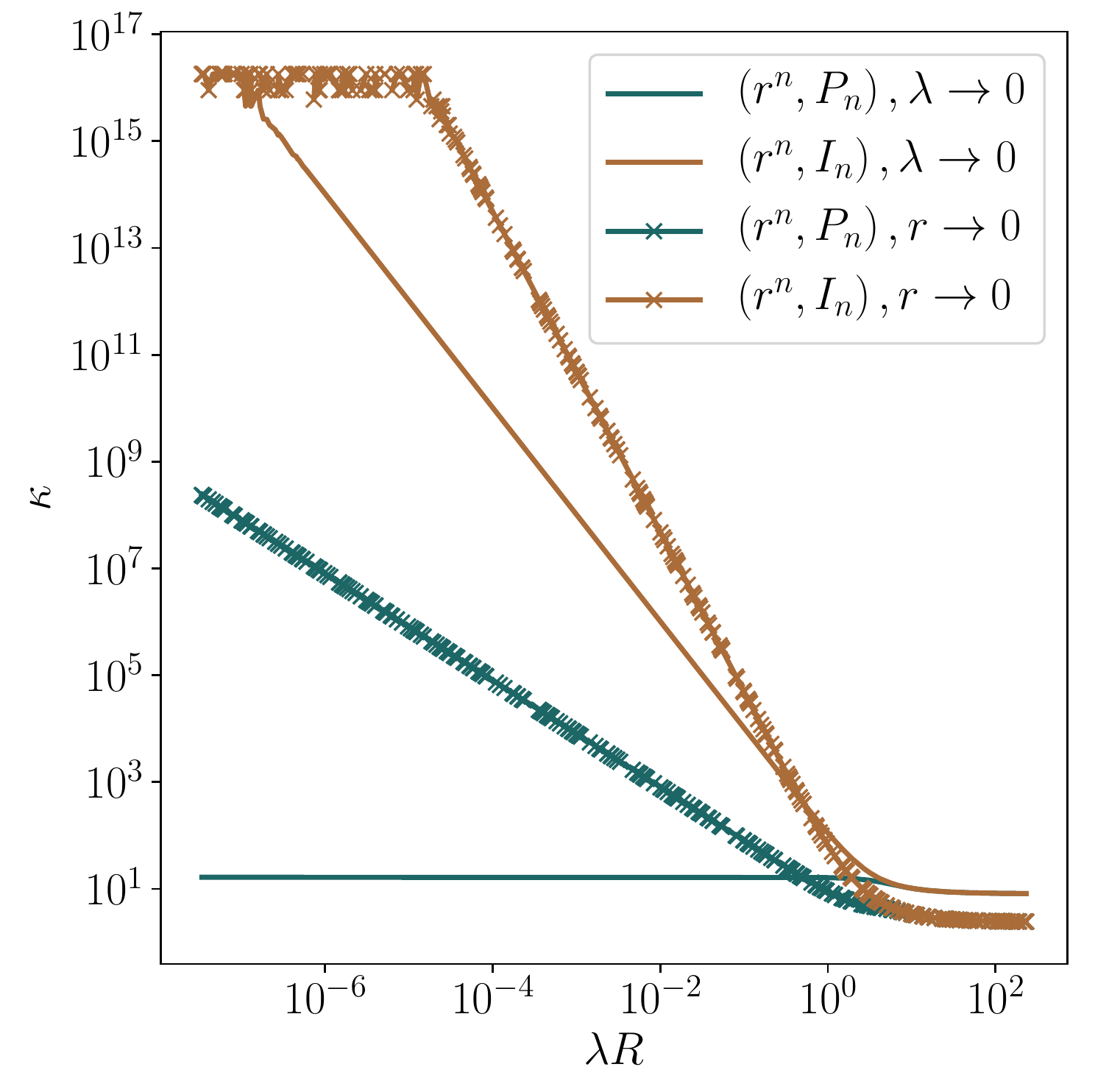}
    \caption{$n=2$}
  \end{subfigure}
  \begin{subfigure}[t]{0.5\textwidth}
    \includegraphics[width=\linewidth]{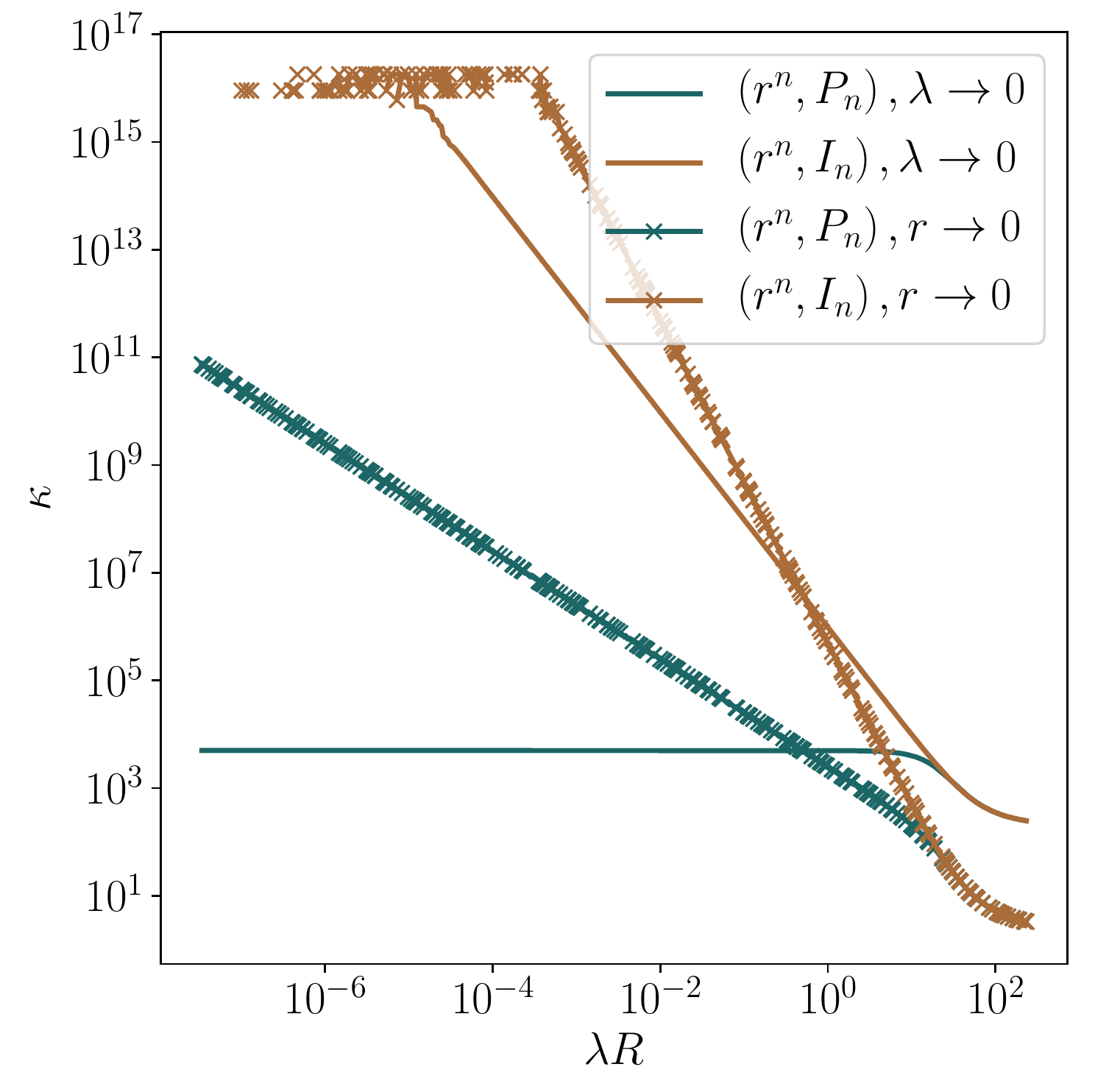}
    \caption{$n=49$}            
  \end{subfigure}

  \caption{We plot the condition number of the linear
    system \cref{eq:matgeneric} (after scaling the columns)
    for various values of $n$ as a function of
    $\lambda R$.}
  \label{fig:condsint}
  
\end{figure}

In \cref{fig:condsint}, we plot the results for the
interior problem. The new basis $(r^{|n|},P_n)$ has smaller
condition numbers than the basis $(r^{|n|},I_n)$, as either
$\lambda$ or $R$ goes to
zero. In the limit as $\lambda$ goes to zero, we see that
the condition number remains roughly constant for
$(r^{|n|},P_n)$. In the limit as $R$ goes to zero, there
is some growth in the condition number for $(r^{|n|},P_n)$,
except for the case $n=0$, where it again remains roughly
constant. For the basis $(r^{|n|},I_n)$, the condition number
grows as either $\lambda$ or $R$ tends to zero and
is larger than for the basis $(r^{|n|},P_n)$. The
condition number tends to grow faster as $R$ goes to
zero compared to the growth as $\lambda$ goes to zero,
except for the case $n = 0$, in which we see the opposite
trend. These results agree well with the analysis of
\cref{sec:analysis}. 

\begin{figure}[h!]
  \begin{subfigure}[t]{0.5\textwidth}
    \includegraphics[width=\linewidth]{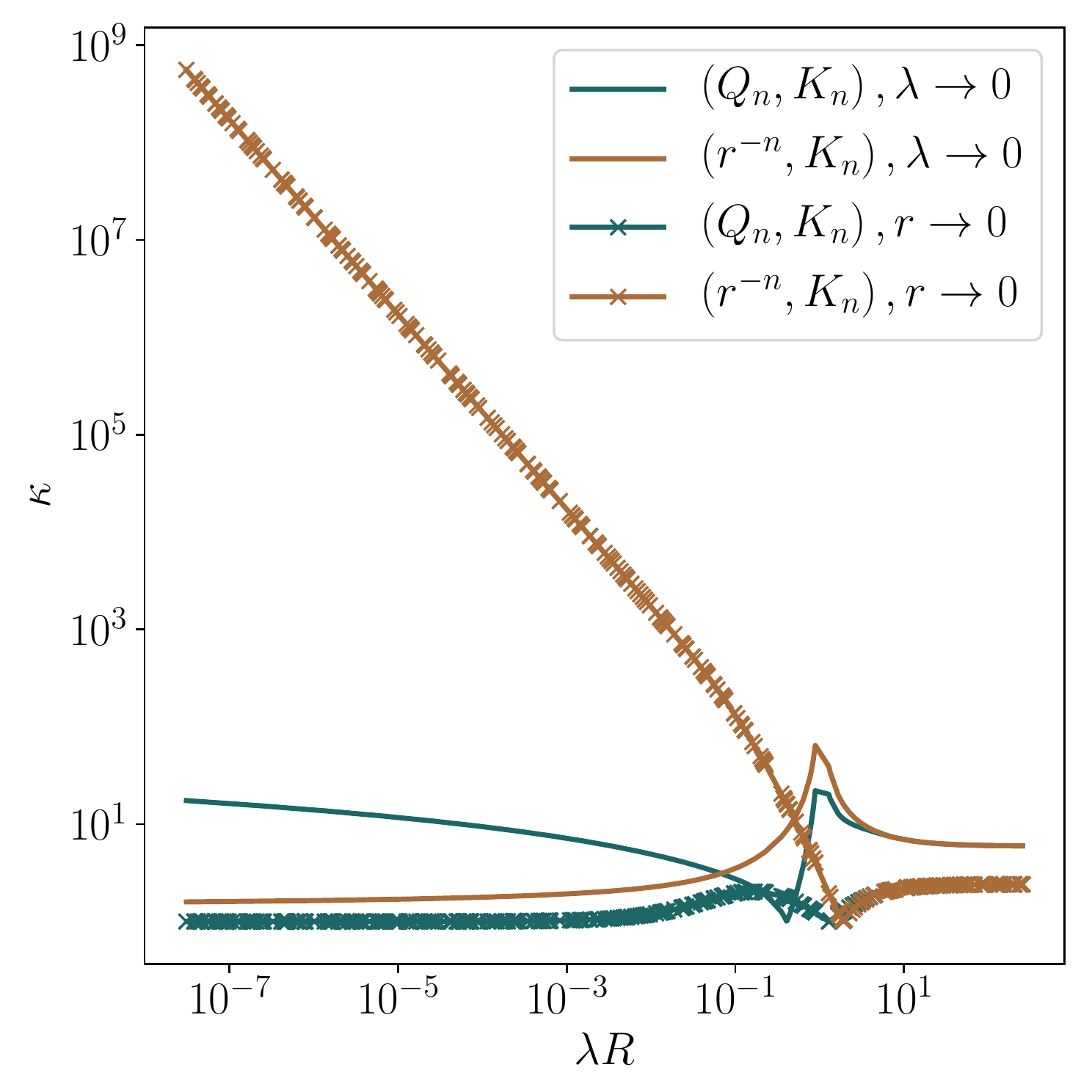}
    \caption{$n = 0$}
  \end{subfigure}
  \begin{subfigure}[t]{0.5\textwidth}
    \includegraphics[width=\linewidth]{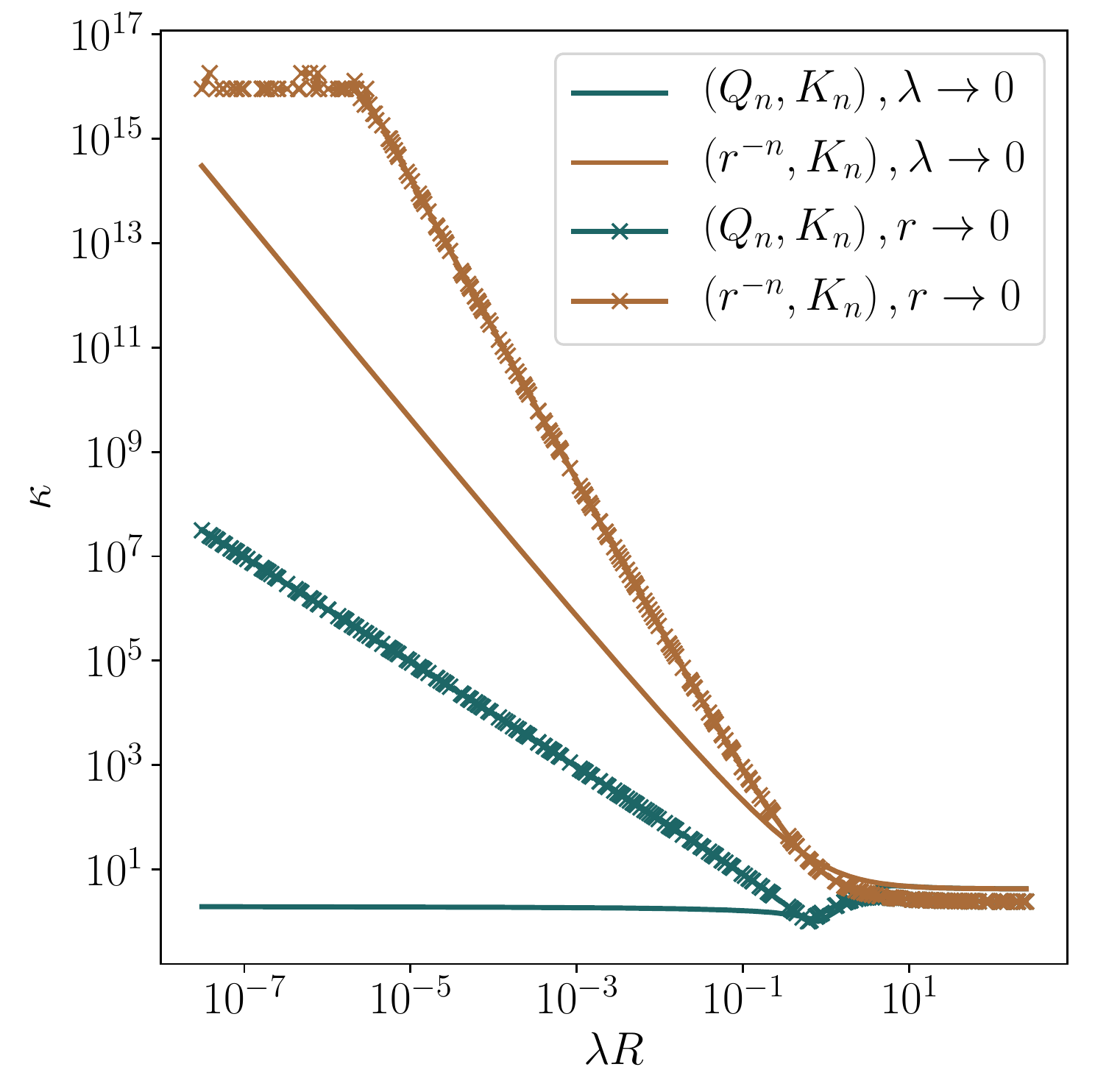}
    \caption{$n=1$}    
  \end{subfigure}
  \begin{subfigure}[t]{0.5\textwidth}
    \includegraphics[width=\linewidth]{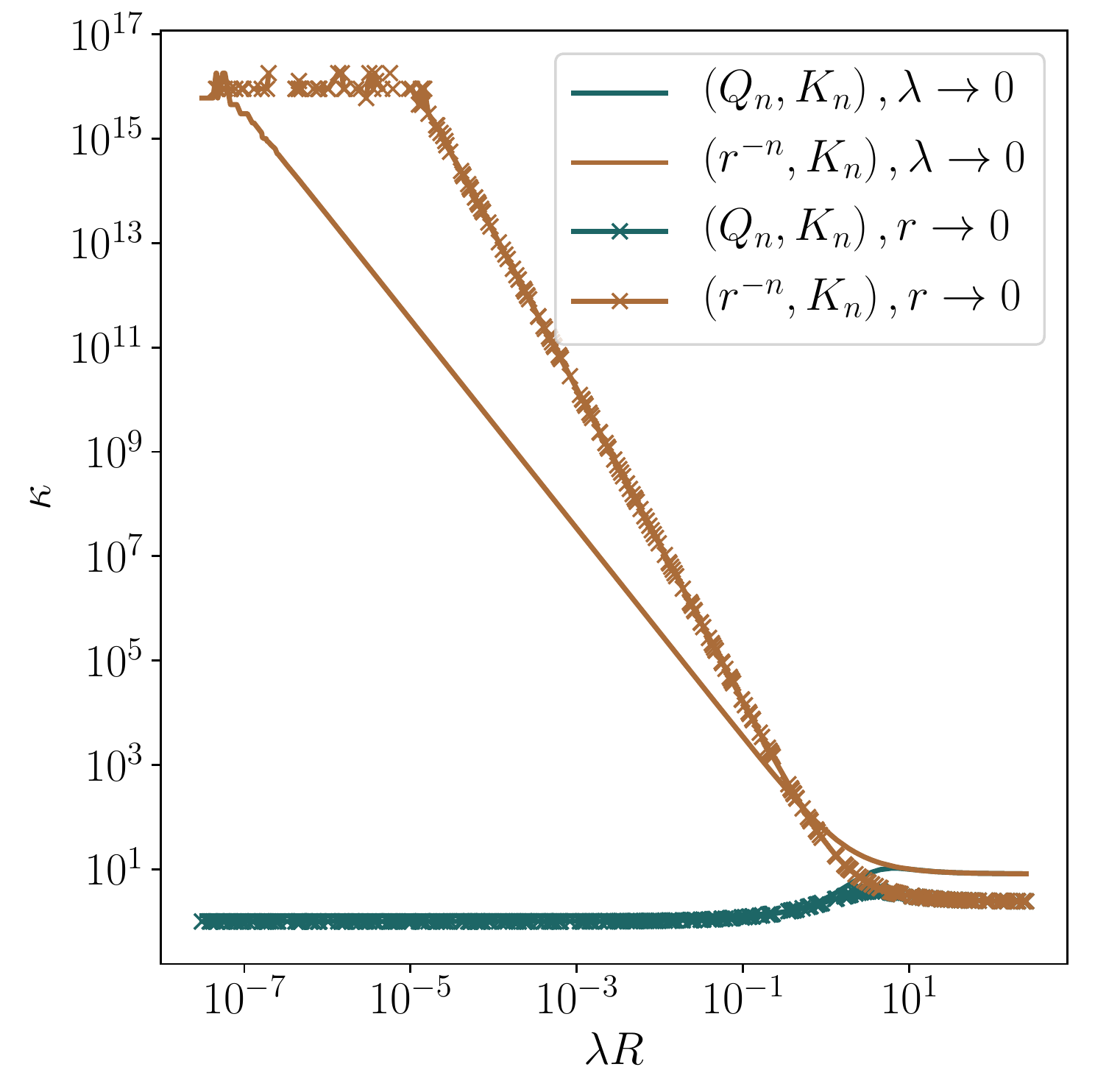}
    \caption{$n=2$}
  \end{subfigure}
  \begin{subfigure}[t]{0.5\textwidth}
    \includegraphics[width=\linewidth]{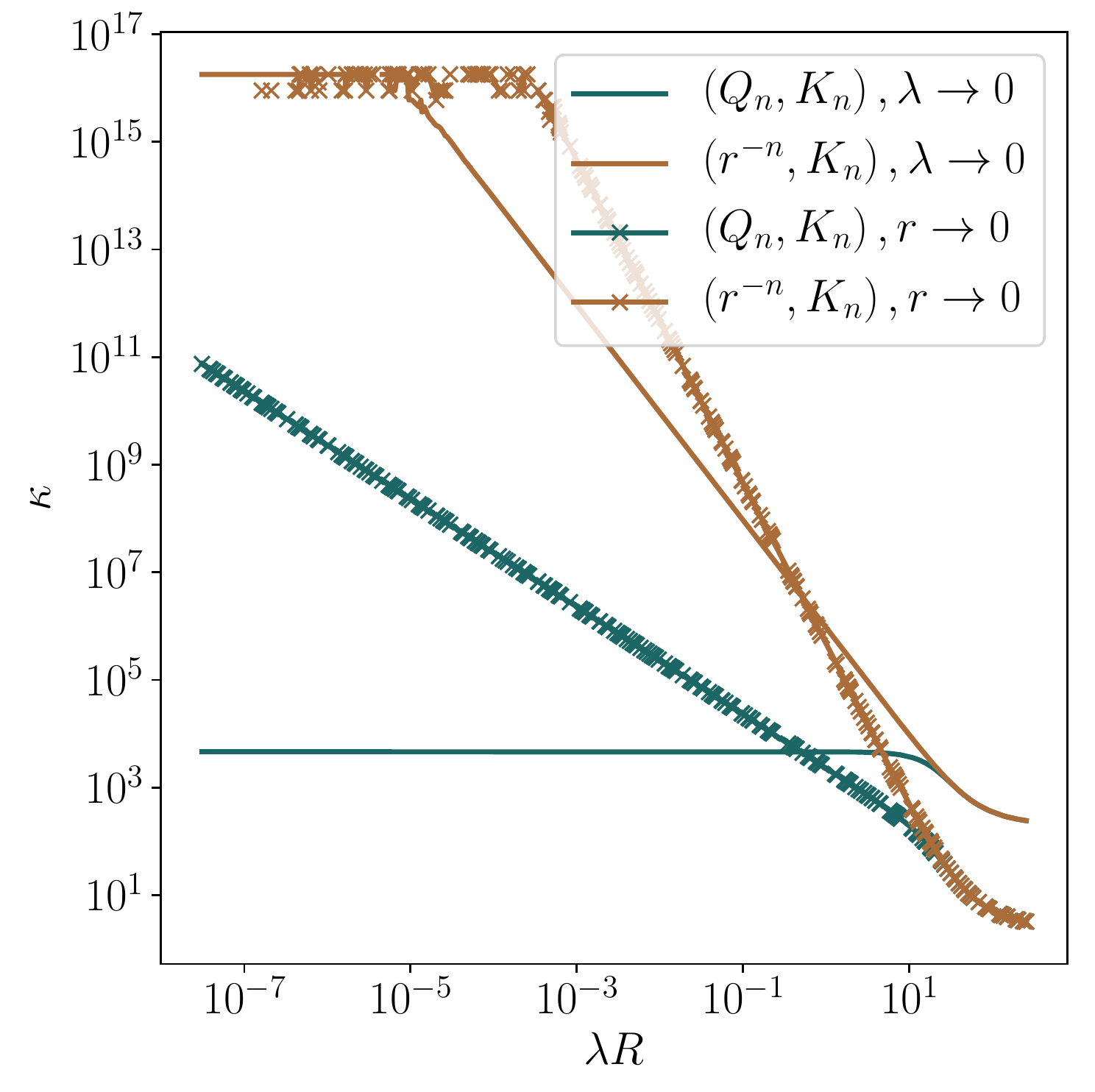}
    \caption{$n=49$}            
  \end{subfigure}

  \caption{We plot the condition number of the linear
    system \cref{eq:matgeneric} (after scaling the columns)
    for various values of $n$ as a function of $\lambda R$.}
  \label{fig:condsext}
  
\end{figure}

In \cref{fig:condsext}, we plot the results for the
exterior problem. The new basis $(Q_n,K_n)$
generally has smaller condition numbers than the
basis $(r^{-|n|},K_n)$, as either $\lambda$ or $R$ goes
to zero.
In the limit as $\lambda$
goes to zero, we see that the condition number remains
roughly constant for $(Q_n,K_n)$. In the limit as $R$
goes to zero, there is some growth in the condition
number for $(Q_n,K_n)$, except for the cases $n=0$
and $n=2$, where it again remains roughly constant.
For the basis $(r^{-|n|},K_n)$, the
condition number grows faster as $R$ goes to
zero compared to the growth as $\lambda$ goes to zero.
The condition number for the basis $(r^{-|n|},K_n)$
does grow as $\lambda$ tends
to zero, except in the case that $n=0$, for which
the condition number remains roughly constant.
Note that for $n=0$ we have used the basis
$(\log(r),K_0)$ as this set of functions is required
to represent the synthetic solution used for the
error analysis of the next section.
Again, these results agree well with the analysis of
\cref{sec:analysis}.

\begin{remark}
  In \cref{fig:condsint,fig:condsext}, we observe
  two distinct behaviors when $\lambda$ tends to
  zero for a fixed $R$ and vice-versa. This phenomenon
  is related to the two relevant scales of the
  bi-Helmholtz kernel noted in \cref{rmk:bihelm},
  though more extreme: whereas the bi-Helmholtz
  operator is the composition of two operators
  with their own dimensionless parameters, the
  modified biharmonic operator is the composition
  of two operators, one of which is scale-invariant
  and the other is not.
  In the next set of experiments, we see that,
  in many applications, there is only one relevant
  scaling for the modified biharmonic equation.
\end{remark}

\subsection{Error plots for known solution}

In order to test the practical effect of the
ill-conditioning seen in the last section,
we use the separation of variables procedure
to solve the modified biharmonic equation
\cref{eq:modbh} on a disk with boundary conditions
corresponding to a known solution. We construct 
this solution using the free-space
Green's function for the modified biharmonic
equation, which is defined as

\begin{equation}
  \mathcal{G}(\x,\y) = -\dfrac1{2\pi \lambda^2} \left ( K_0(\lambda \rho) +
  \log(\rho) \right ) \; ,
\end{equation}
where $\rho = \sqrt { (x_1-y_1)^2 + (x_2-y_2)^2 }$.
The solution is then set to be

\begin{equation}
  u(\x; \lambda) = \sum_{j=1}^{n_s} \lambda^2 c_j \mathcal{G}(\x,\s_j)
  + \lambda d_j \partial_{v_{j,1}} \mathcal{G}(\x,\s_j) +
  q_j \partial_{v_{j,2}v_{j,3}} \mathcal{G}(\x,\s_j) \, ,
\end{equation}
where the $\s_j$ are $n_s = 100$ ``source'' points
located outside
of the domain, the $c_j$ are drawn uniformly randomly
from $[-1,1]$, the $d_j$ and $q_j$ are drawn uniformly
randomly from $[0,1]$, and the vectors $v_{j,1},v_{j,2},v_{j,3}$
are defined by drawing the entries uniformly at random
from $[-1/2,1/2]$ and normalizing. We note that the
$\lambda^2$ and $\lambda$ scales in front of the $c_j$
and $d_j$ are included to ensure that these
terms are of roughly the same size as $\lambda$ shrinks.

To implement separation of variables, we discretize the
boundary with $M = 100$ points, so that $N=49$ (the 
separation of variables expansion runs from 
$-N$ to $N+1$ as in \cref{sec:sep}). We
evaluate the function $u(\x; \lambda)$ and its normal
derivative on the boundary
of the disk and compute their Fourier coefficients, i.e.
the values $f_n$ and $g_n$ as in \cref{eq:fn,eq:gn},
using the FFT. We then solve for the coefficients
$\alpha_n$ and $\beta_n$ by inverting the linear system

\begin{equation}
\begin{pmatrix}
    F_n(R) & G_n(R) \\
    F_n'(R) & G_n'(R)
\end{pmatrix}
\begin{pmatrix} \alpha_n \\ \beta_n 
\end{pmatrix}
=\begin{pmatrix}
f_n \\ g_n
\end{pmatrix} \; ,
\end{equation}
where $(F_n(r),G_n(r))$ is an appropriate pair of basis
functions. For the sake of stability, the inversion
is performed using Gaussian elimination with complete
pivoting. Once the coefficients $\alpha_n$ and $\beta_n$
have been computed, the approximate solution $\hat{u}$
can be evaluated using the formula

\begin{equation}
  \hat{u}(\x) = \sum_{n=-N}^{N+1} (\alpha_n F_n(\rho) +
  \beta_n G_n(\rho) ) e^{in\theta} \; ,
\end{equation}
where $(\rho,\theta)$ are the polar coordinates of
$\x$. The derivatives of $\hat{u}$ can be obtained
by differentiating this expression.

To measure the performance of the separation of
variables method, we evaluate $\hat{u}$ and its
first and second derivatives at $n_t = 100$
``target'' points $\t_i$ located inside the domain.
We then define three error measures

\begin{align}
  E_u &= \sqrt{\dfrac{\sum_{i=1}^{n_t} (u(\t_i ; \lambda)-\hat{u}(\t_i))^2 }
  {\sum_{i=1}^{n_t} u(\t_i ; \lambda)^2 }} \; ,\\
  E_g &= \sqrt{\dfrac{ \sum_{i=1}^{n_t}
        (u_x(\t_i ; \lambda)-\hat{u}_x(\t_i))^2
      + (u_y(\t_i ; \lambda)-\hat{u}_y(\t_i))^2 }
    { \sum_{i=1}^{n_t} u_x(\t_i ; \lambda)^2
      + u_y(\t_i ; \lambda)^2 }} \; , \\
  E_h &= \sqrt{\dfrac{ \sum_{i=1}^{n_t}
        (u_{xx}(\t_i ; \lambda)-\hat{u}_{xx}(\t_i))^2
      + (u_{xy}(\t_i ; \lambda)-\hat{u}_{xy}(\t_i))^2
      + (u_{yy}(\t_i ; \lambda)-\hat{u}_{yy}(\t_i))^2      
    }
    { \sum_{i=1}^{n_t} u_{xx}(\t_i ; \lambda)^2
      + u_{xy}(\t_i ; \lambda)^2
      + u_{yy}(\t_i ; \lambda)^2 } } \; ,
\end{align}
which represent the relative error in the solution,
gradient, and Hessian, respectively.

As in the previous
section, we run the tests with a variety of values for
the radius of the disk $R$ and the parameter $\lambda$.
To observe the effect of
changing the radius of the domain, we run experiments
with $\lambda$ fixed (at $\lambda = 0.5$) and, for
each $j$ from $-24$ to $8$, ten values of $R$ drawn
uniformly at random from the interval $[2^{j},2^{j+1}]$.
Likewise, to observe the effect of changing the
parameter $\lambda$, we run experiments with $R$ fixed
(at $R = 0.5$) and, for each $j$ from $-24$ to $8$,
ten values of $\lambda$ drawn uniformly at random from
the interval $[2^{j},2^{j+1}]$. In order to compare these,
we plot the error measures $E_u$, $E_g$, and $E_h$
as functions of the product $\lambda R$.

\begin{figure}[h!]

  \begin{subfigure}[t]{0.5\textwidth}
    \centering
    \includegraphics[width=\linewidth]{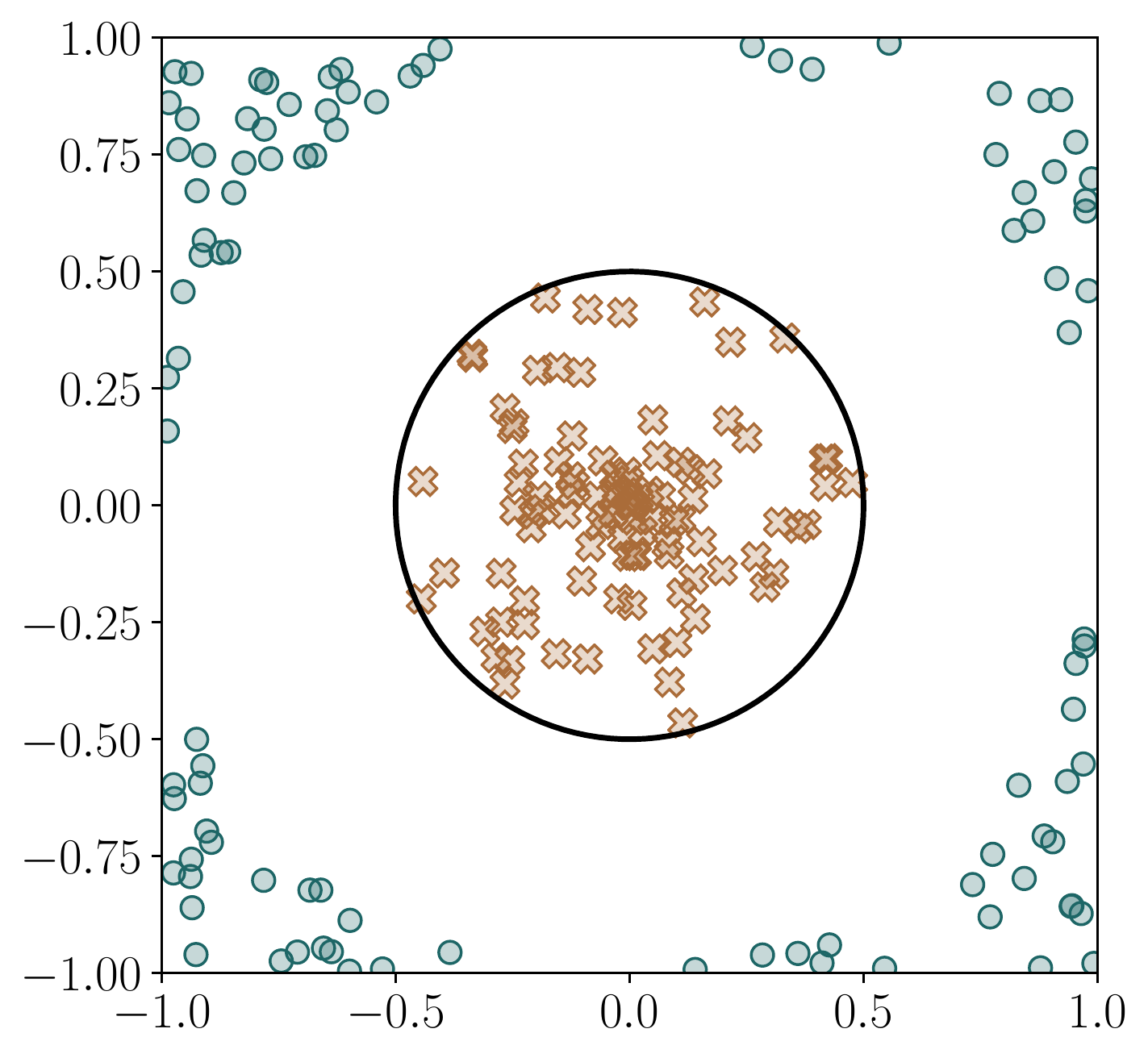}
    \caption{Interior problem geometry}
    \label{sfig:intgeo}
  \end{subfigure}
  \begin{subfigure}[t]{0.5\textwidth}
    \centering
    \includegraphics[width=\linewidth]{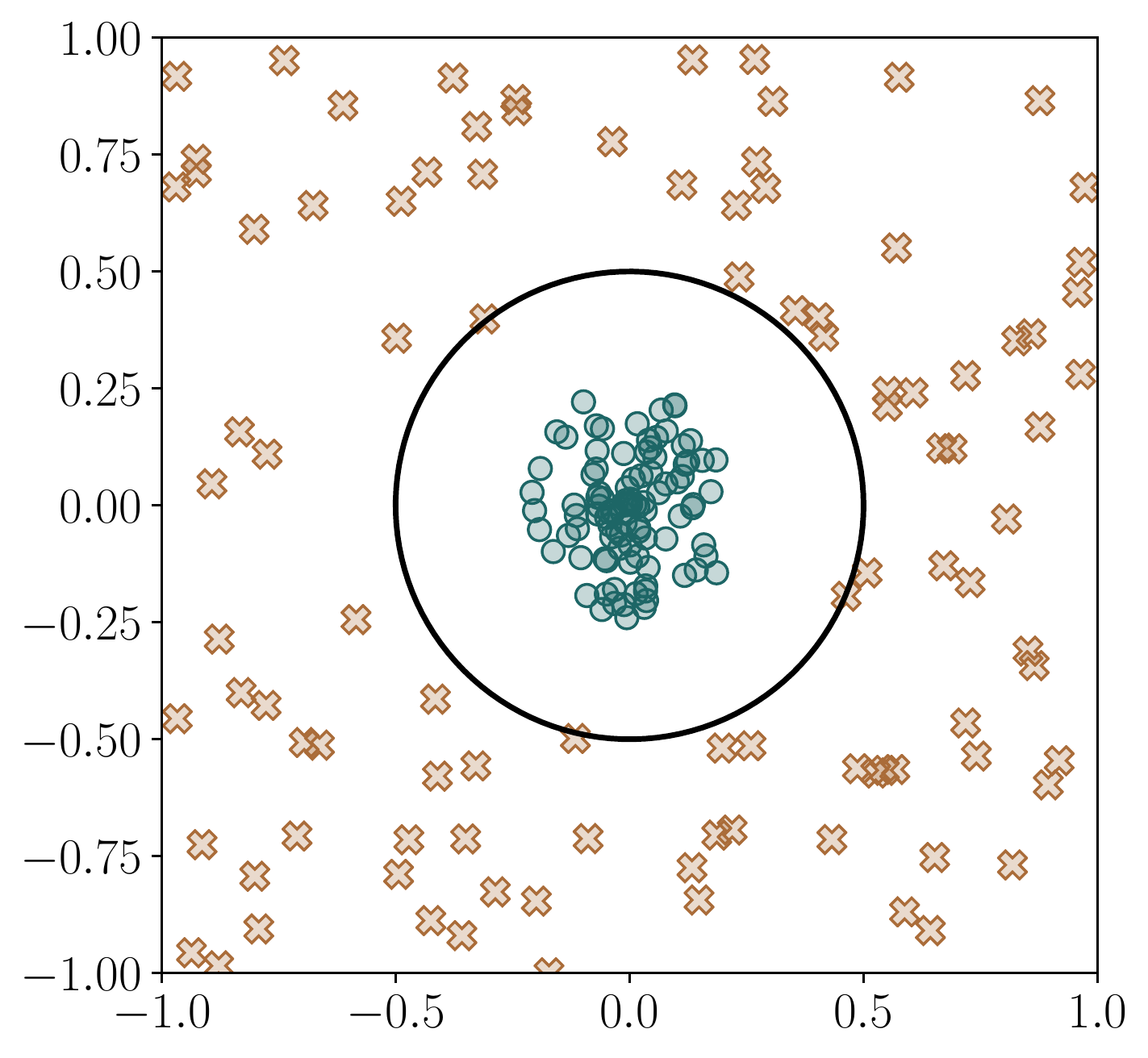}
    \caption{Exterior problem geometry}
    \label{sfig:extgeo}    
  \end{subfigure}
  \caption{Sample geometries for the interior
    and exterior problems with $R=0.5$.
    The targets are marked by crosses and the sources
    by circles.}

\end{figure}

For the interior problem, the source points are drawn
uniformly at random from the box $[-2R,2R] \times [-2R,2R]$
outside of the disk of radius $2R$. With this placement
of source points, the length $N=49$ expansion for
$\hat{u}$ should be sufficient for approximately machine
precision accuracy, based on standard multipole estimates
\cite{greengard1987fast,cheng2006adaptive}. The
target points are drawn uniformly at random from the
disk of radius $R$. We plot a sample arrangement of
the source and target points for the interior problem
in \cref{sfig:intgeo}.

Similarly, for the exterior problem, the source
points are drawn uniformly at random from the
disk of radius $R/2$. This placement is again chosen
so that the length $N= 49$ expansion for $\hat{u}$
is sufficient for approximately machine precision
accuracy. The target points are drawn uniformly
at random from the box $[-2R,2R] \times [-2R,2R]$
outside of the disk of radius $R$. We plot a sample
arrangement of the source and target points for the
exterior problem in \cref{sfig:extgeo}.

\begin{figure}[h!]

  \begin{subfigure}[t]{0.32\textwidth}
    \includegraphics[width=\linewidth]{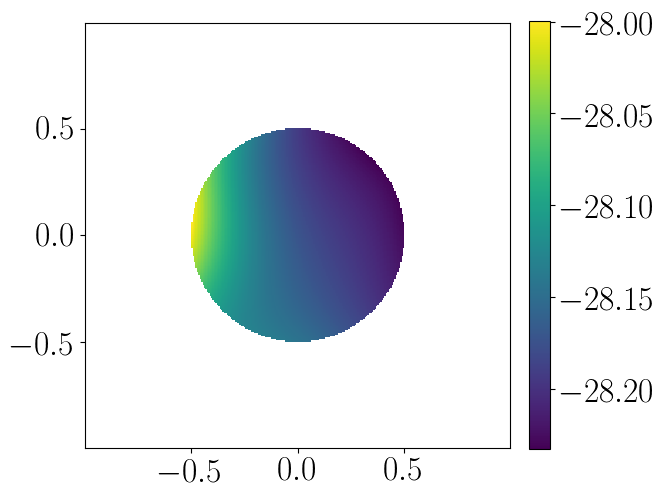}
    \caption{$u$}
  \end{subfigure}
  \begin{subfigure}[t]{0.32\textwidth}
    \includegraphics[width=\linewidth]{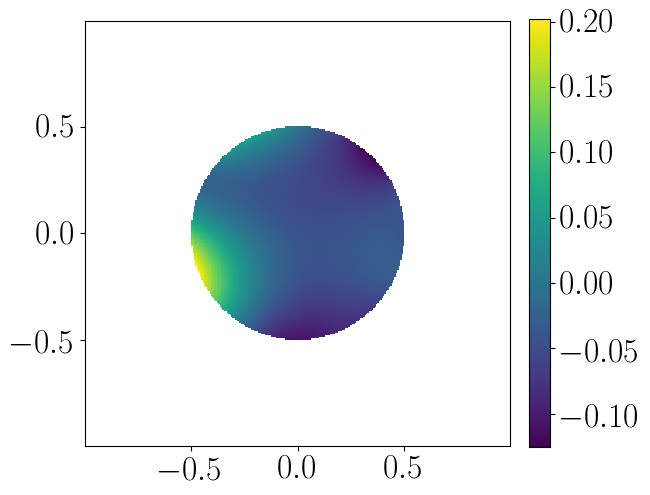}
    \caption{$u_x$}    
  \end{subfigure}
  \begin{subfigure}[t]{0.32\textwidth}
    \includegraphics[width=\linewidth]{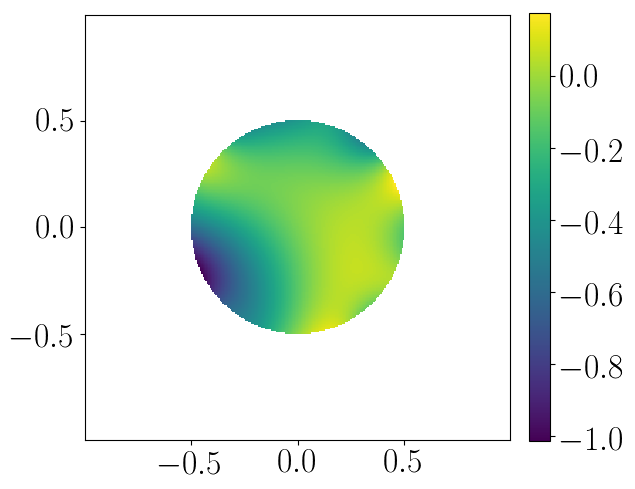}
    \caption{$u_{xy}$}
  \end{subfigure}

  \caption{Heatmaps of the exact solution $u$ and
    select derivatives for the interior problem with
    $R= 0.5$ and $\lambda = 2^{-24}$.}

  \label{fig:intsol}  
  
\end{figure}

\begin{figure}[h!]

  \begin{subfigure}[t]{0.32\textwidth}
    \includegraphics[width=\linewidth]{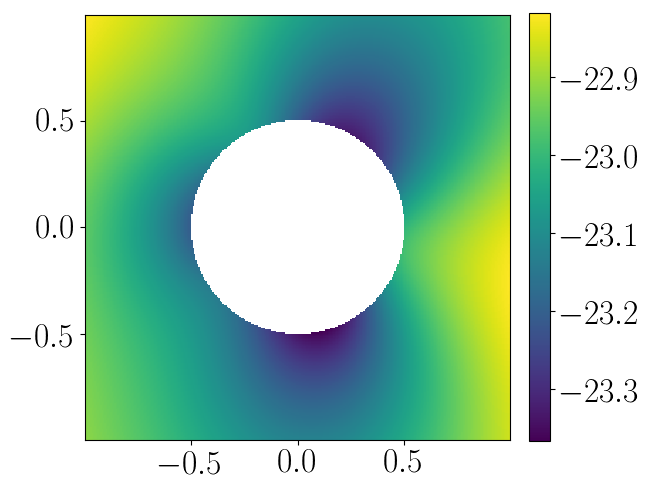}
    \caption{$u$}
  \end{subfigure}
  \begin{subfigure}[t]{0.32\textwidth}
    \includegraphics[width=\linewidth]{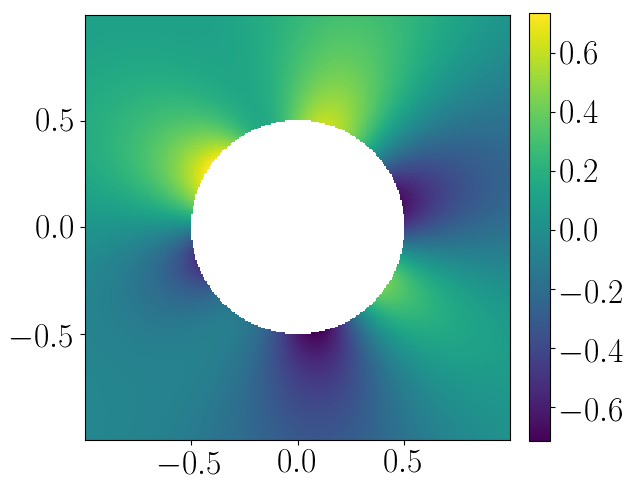}
    \caption{$u_x$}    
  \end{subfigure}
  \begin{subfigure}[t]{0.32\textwidth}
    \includegraphics[width=\linewidth]{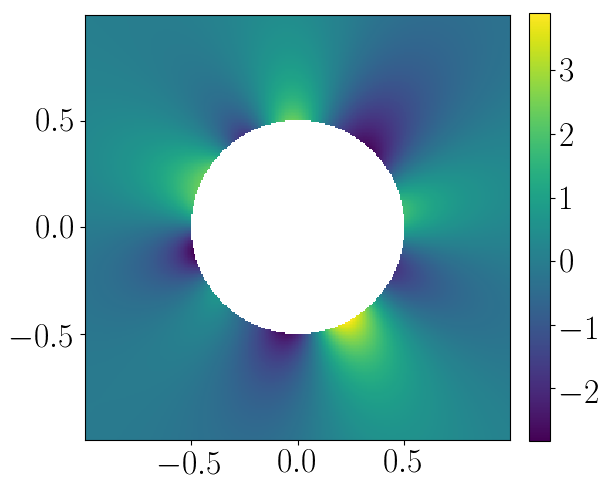}
    \caption{$u_{xy}$}
  \end{subfigure}

  \caption{Heatmaps of the exact solution $u$ and
    select derivatives for the exterior problem with
    $R= 0.5$ and $\lambda = 2^{-24}$.}

  \label{fig:extsol}
  
\end{figure}

In \cref{fig:intsol,fig:extsol} we plot a sample
exact solution $u$ and some select derivatives for
the interior and exterior problems, respectively.
For the following error plots, we consider the
error in approximating these solutions using
the bases we have discussed above as well as using
what we call the ``exact difference''. We include
the exact difference figure to emphasize that
it is not only the ability to recover $\alpha_n$
and $\beta_n$ that causes trouble using the
na\"{i}ve bases. Because $u$ is defined
in terms of the Green's function $\mathcal{G}$, which
is simply a scaled sum of the Green's functions
for the Laplace and modified Helmholtz equations,
we could reasonably evaluate $u$ by evaluating
these parts separately and combining them in the
end. Let $u_L$ and $u_H$ be defined as 

\begin{equation}
  u_L(\x; \lambda) = \dfrac{1}{2\pi \lambda^2}
  \sum_{j=1}^{n_s} \lambda^2 c_j \log \|\x-\s_j\|_2
  + \lambda d_j \partial_{v_{j,1}} \log \|\x -\s_j\|_2 +
  q_j \partial_{v_{j,2}v_{j,3}} \log \| \x - \s_j\|_2
\end{equation}
and 
\begin{equation}
  u_H(\x; \lambda) = -\dfrac{1}{2\pi \lambda^2}
  \sum_{j=1}^{n_s} \lambda^2 c_j K_0 (\lambda \|\x-\s_j\|_2)
  + \lambda d_j \partial_{v_{j,1}} K_0(\lambda \|\x -\s_j\|_2) +
  q_j \partial_{v_{j,2}v_{j,3}} K_0(\lambda \| \x - \s_j\|_2) \, .
\end{equation}
The ``exact difference'' error below is the
error in evaluating $u$ as the difference
$u_H - u_L$ in floating point arithmetic.

\begin{figure}[h!]

  \begin{subfigure}[t]{0.32\textwidth}
    \includegraphics[width=\linewidth]{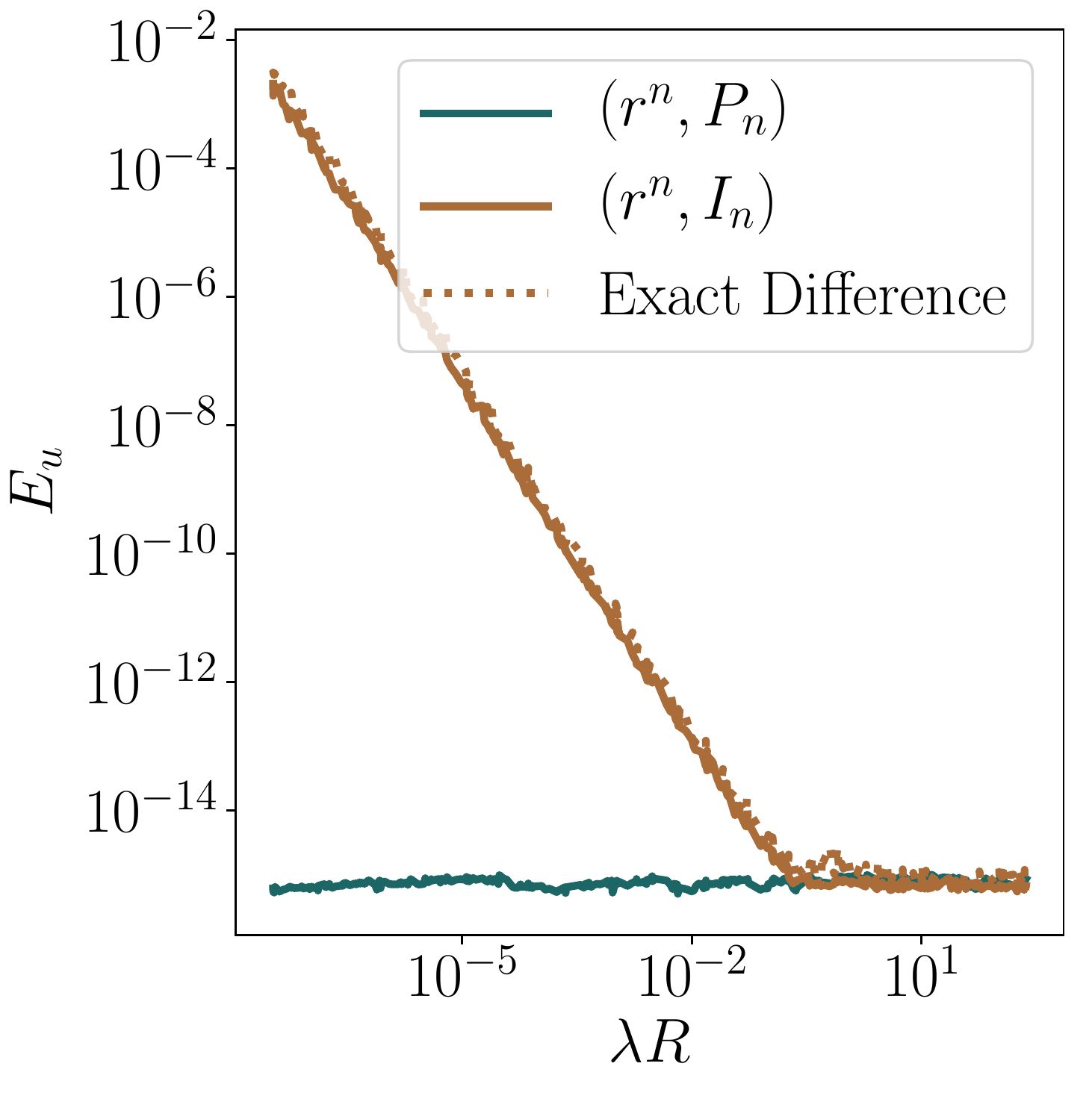}
    \caption{$E_u$ as $\lambda$ goes to zero.}
  \end{subfigure}
  \begin{subfigure}[t]{0.32\textwidth}
    \includegraphics[width=\linewidth]{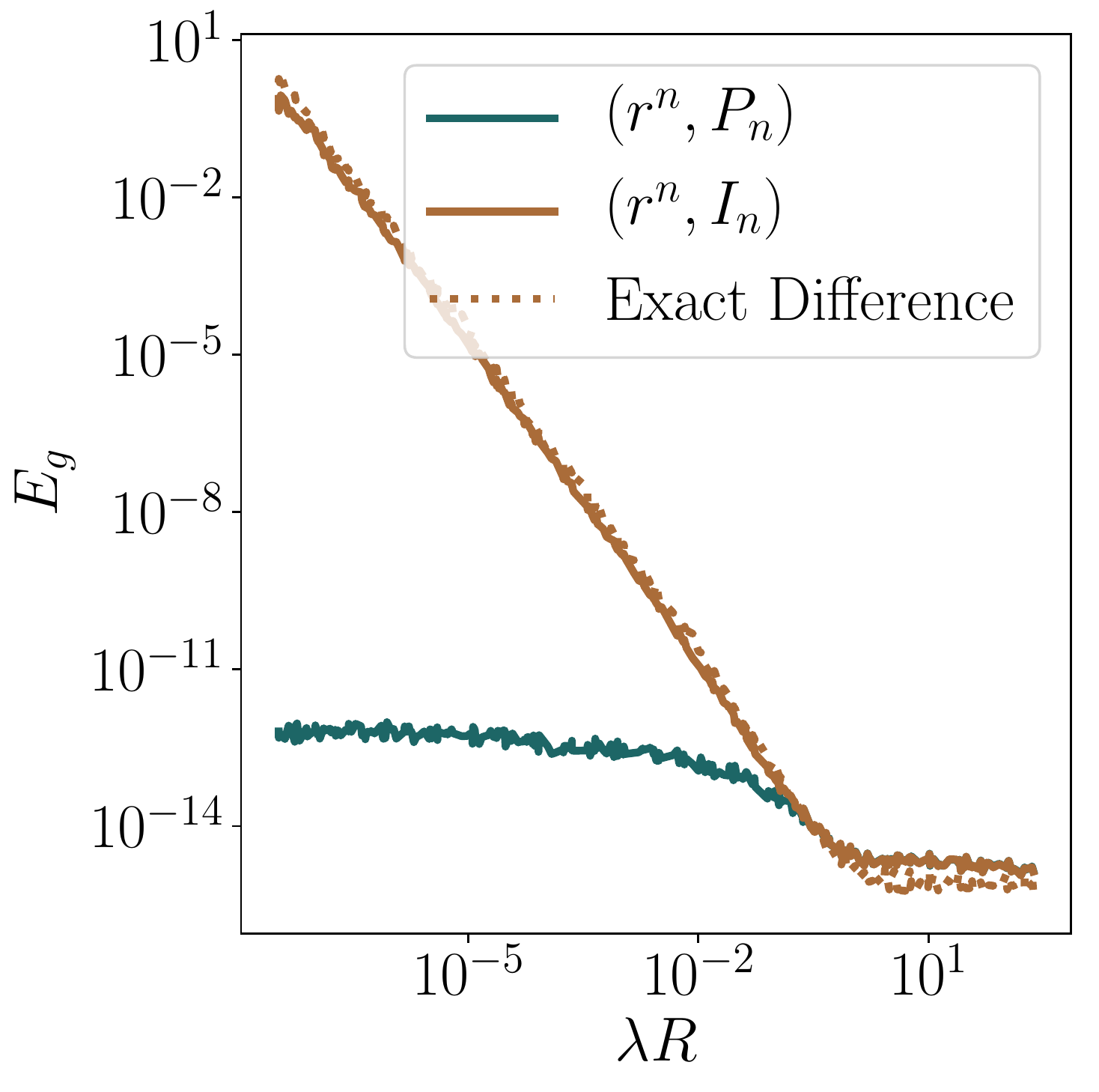}
    \caption{$E_g$ as $\lambda$ goes to zero.}    
  \end{subfigure}
  \begin{subfigure}[t]{0.32\textwidth}
    \includegraphics[width=\linewidth]{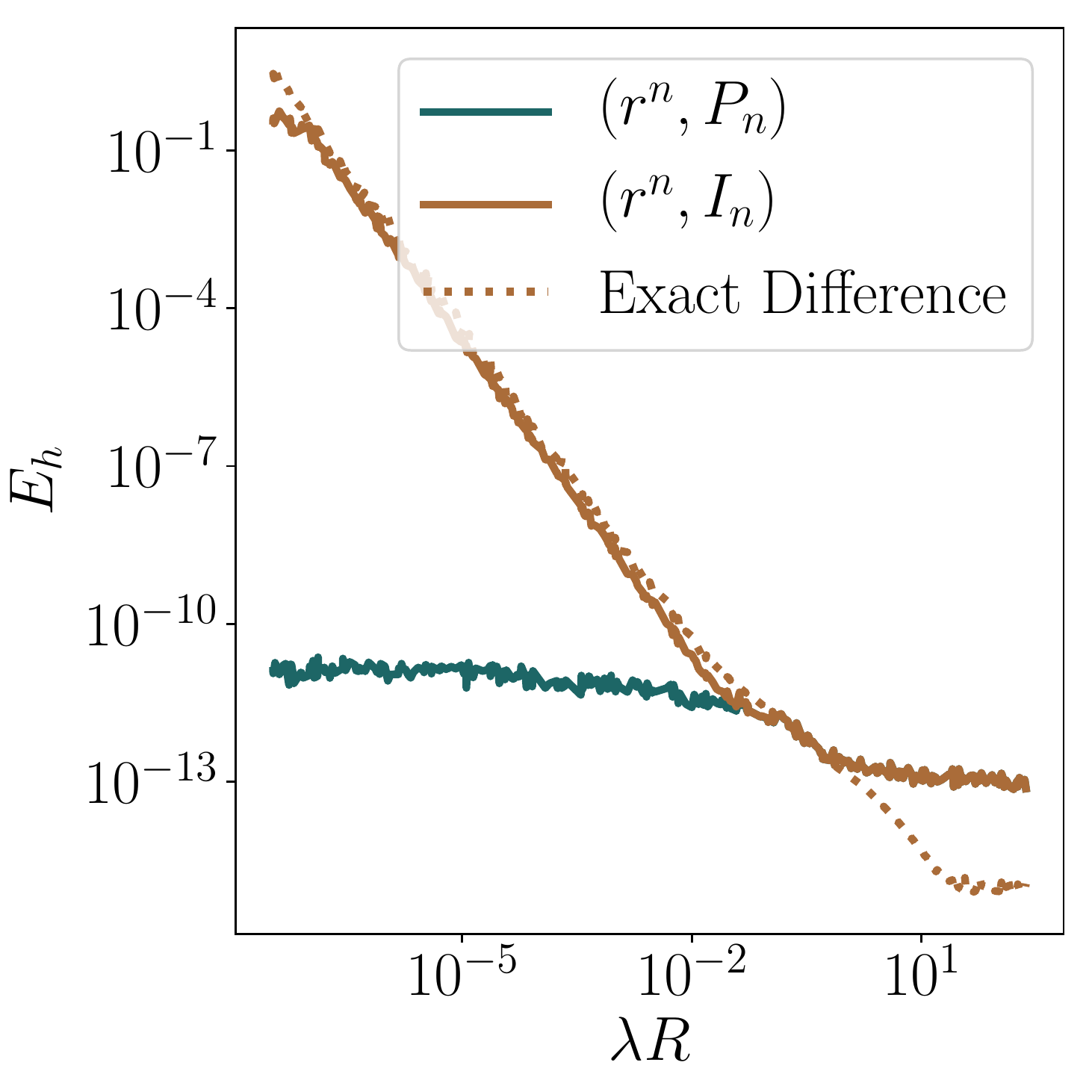}
    \caption{$E_h$ as $\lambda$ goes to zero.}
  \end{subfigure}

  \begin{subfigure}[t]{0.32\textwidth}
    \includegraphics[width=\linewidth]{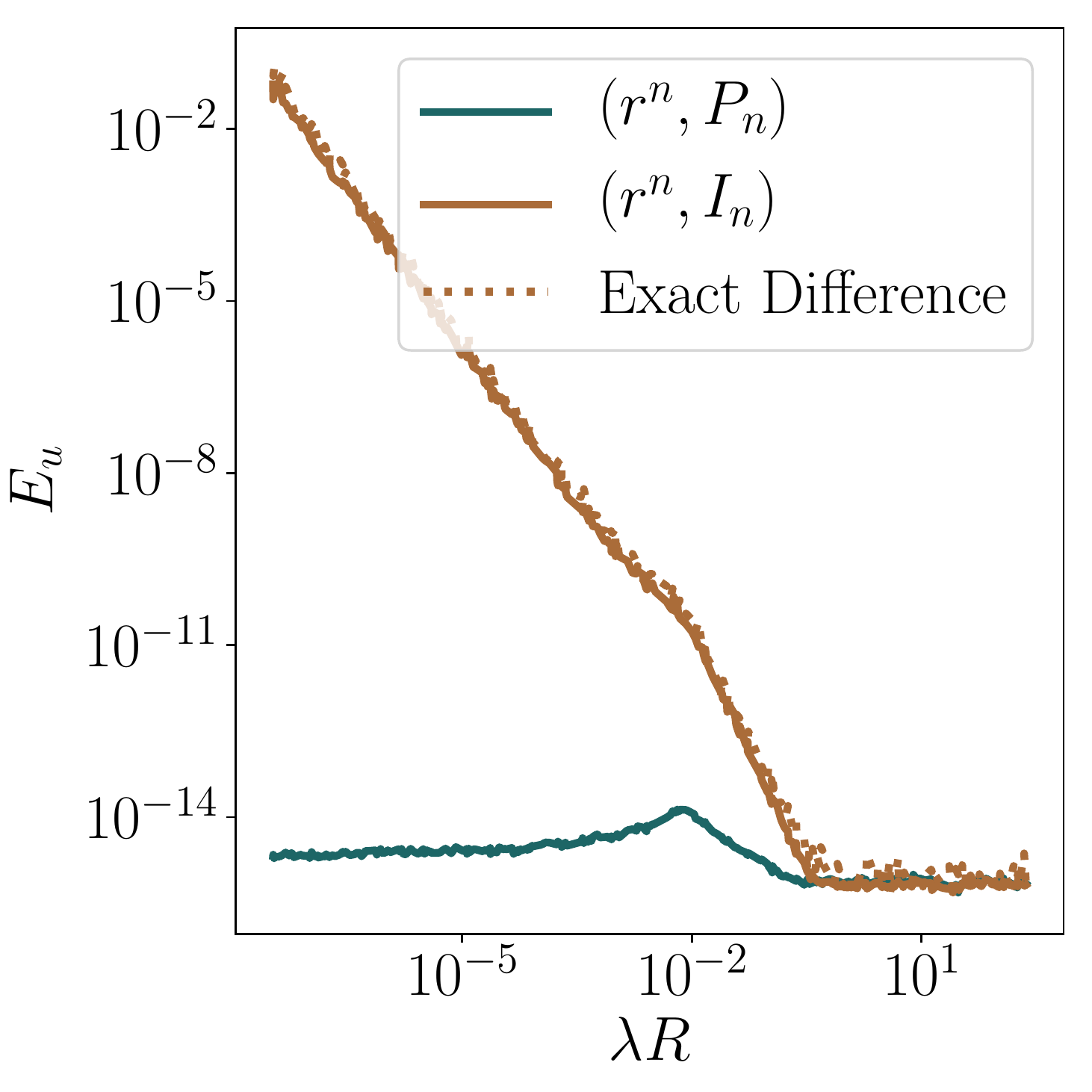}
    \caption{$E_u$ as $R$ goes to zero.}
  \end{subfigure}
  \begin{subfigure}[t]{0.32\textwidth}
    \includegraphics[width=\linewidth]{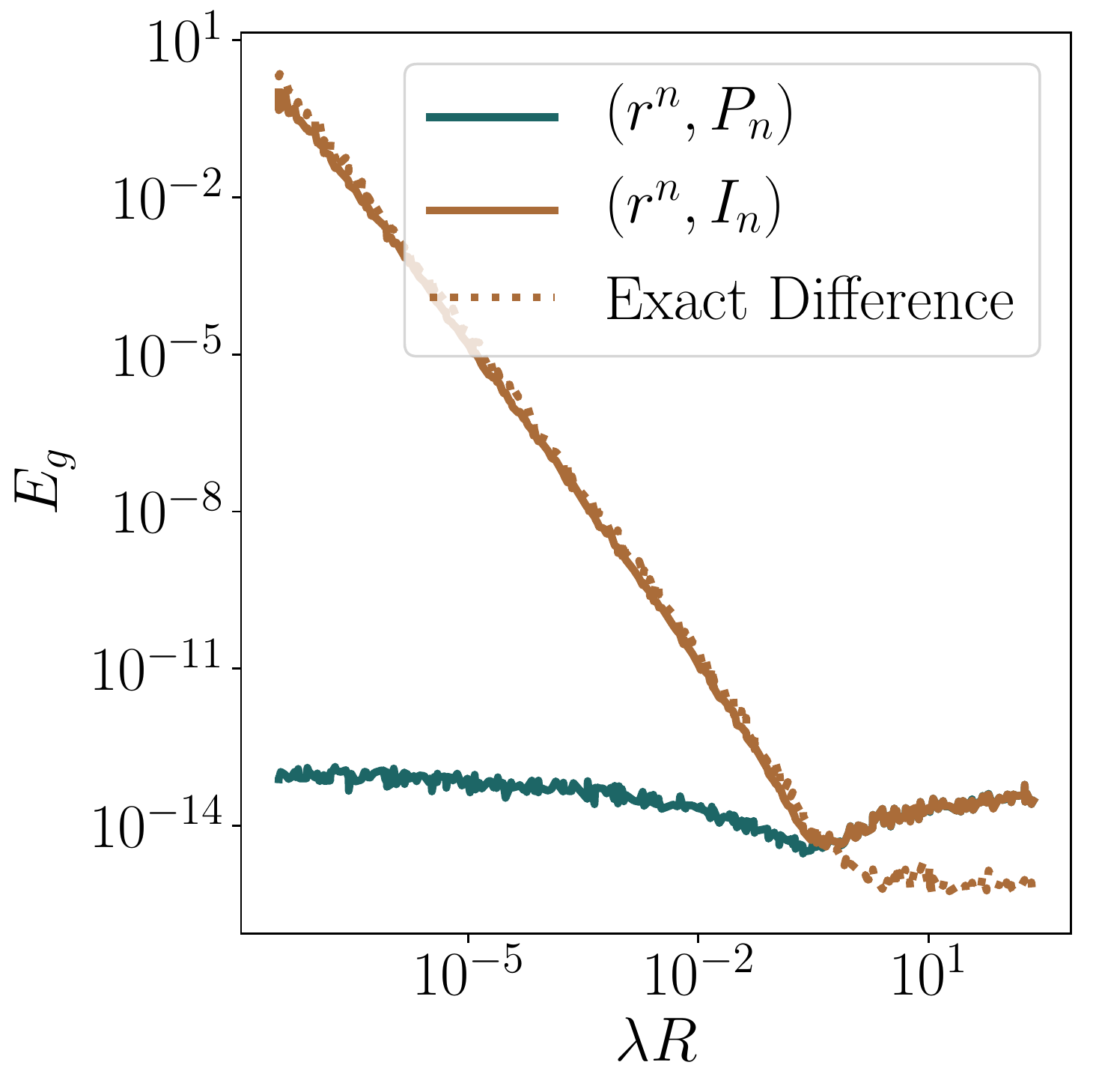}
    \caption{$E_g$ as $R$ goes to zero.}    
  \end{subfigure}
  \begin{subfigure}[t]{0.32\textwidth}
    \includegraphics[width=\linewidth]{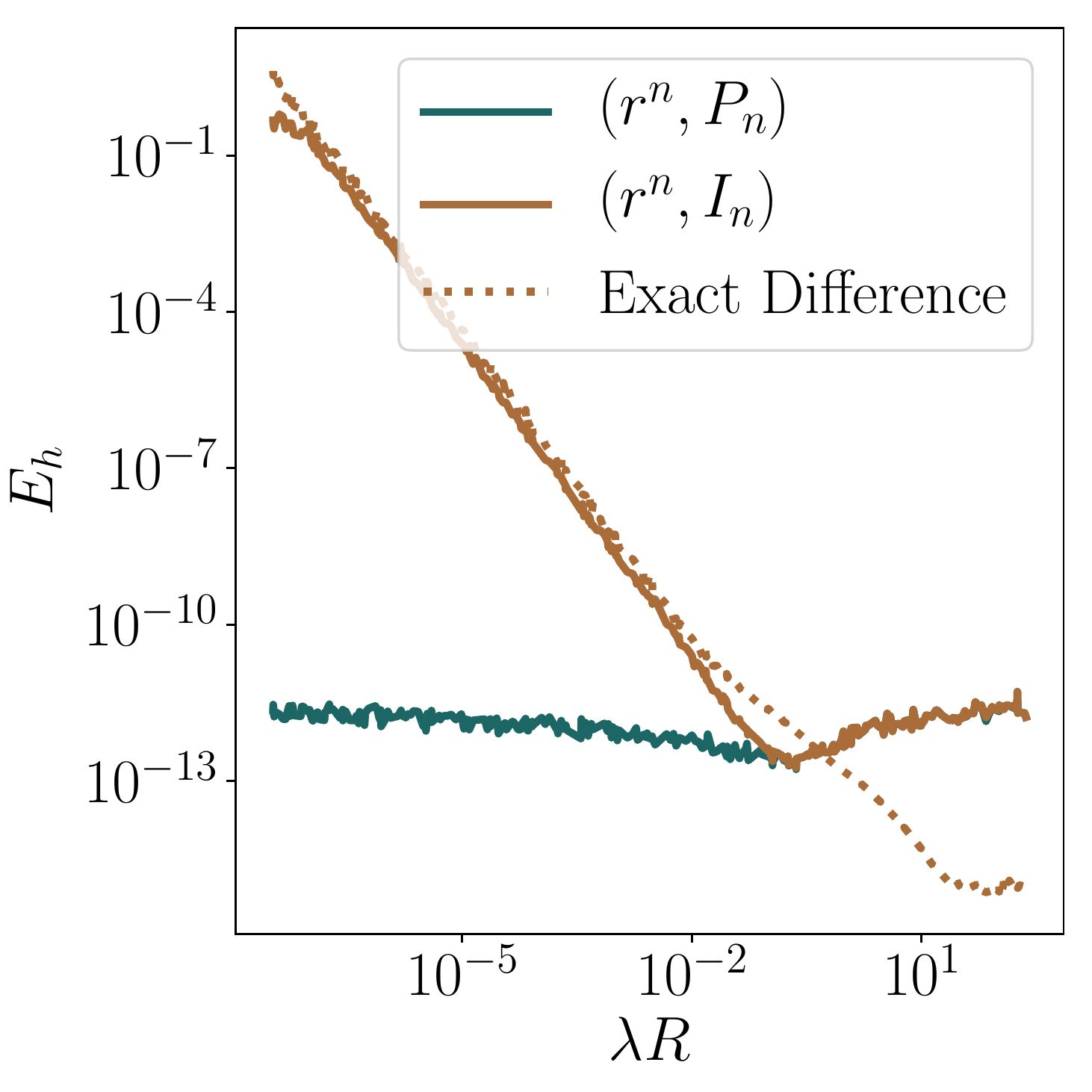}
    \caption{$E_h$ as $R$ goes to zero.}
  \end{subfigure}

  \caption{Interior problem. In the top row, we plot
    the error measures as functions of $\lambda R$
    for $R=0.5$ as $\lambda$
    goes to zero. In the bottom row, we plot the error
    measures as functions of $\lambda R$ for $\lambda = 0.5$
    as $R$ goes to zero.}

  \label{fig:errint}
  
\end{figure}

\begin{figure}[h!]

  \begin{subfigure}[t]{0.32\textwidth}
    \includegraphics[width=\linewidth]{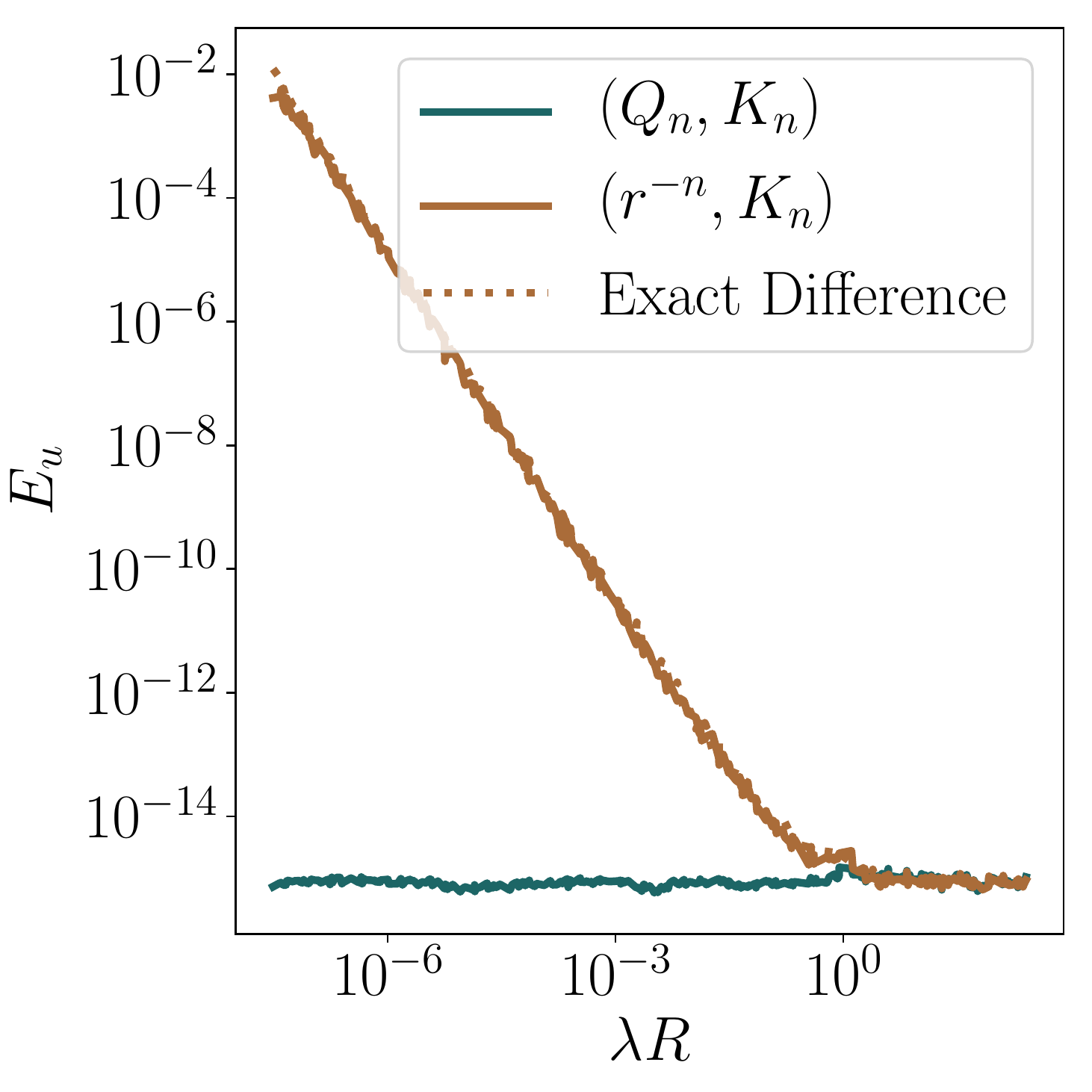}
    \caption{$E_u$ as $\lambda$ goes to zero.}
  \end{subfigure}
  \begin{subfigure}[t]{0.32\textwidth}
    \includegraphics[width=\linewidth]{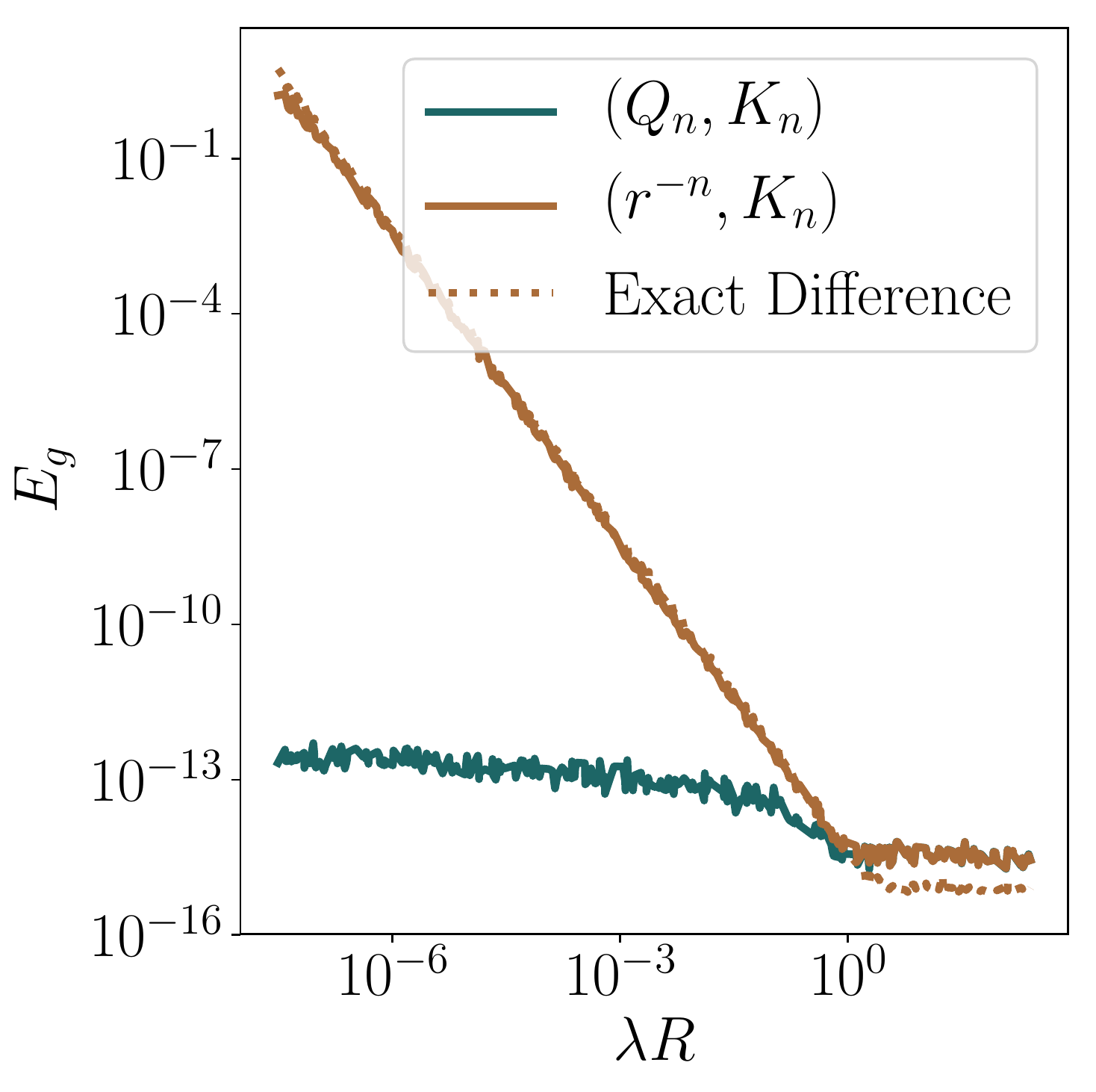}
    \caption{$E_g$ as $\lambda$ goes to zero.}    
  \end{subfigure}
  \begin{subfigure}[t]{0.32\textwidth}
    \includegraphics[width=\linewidth]{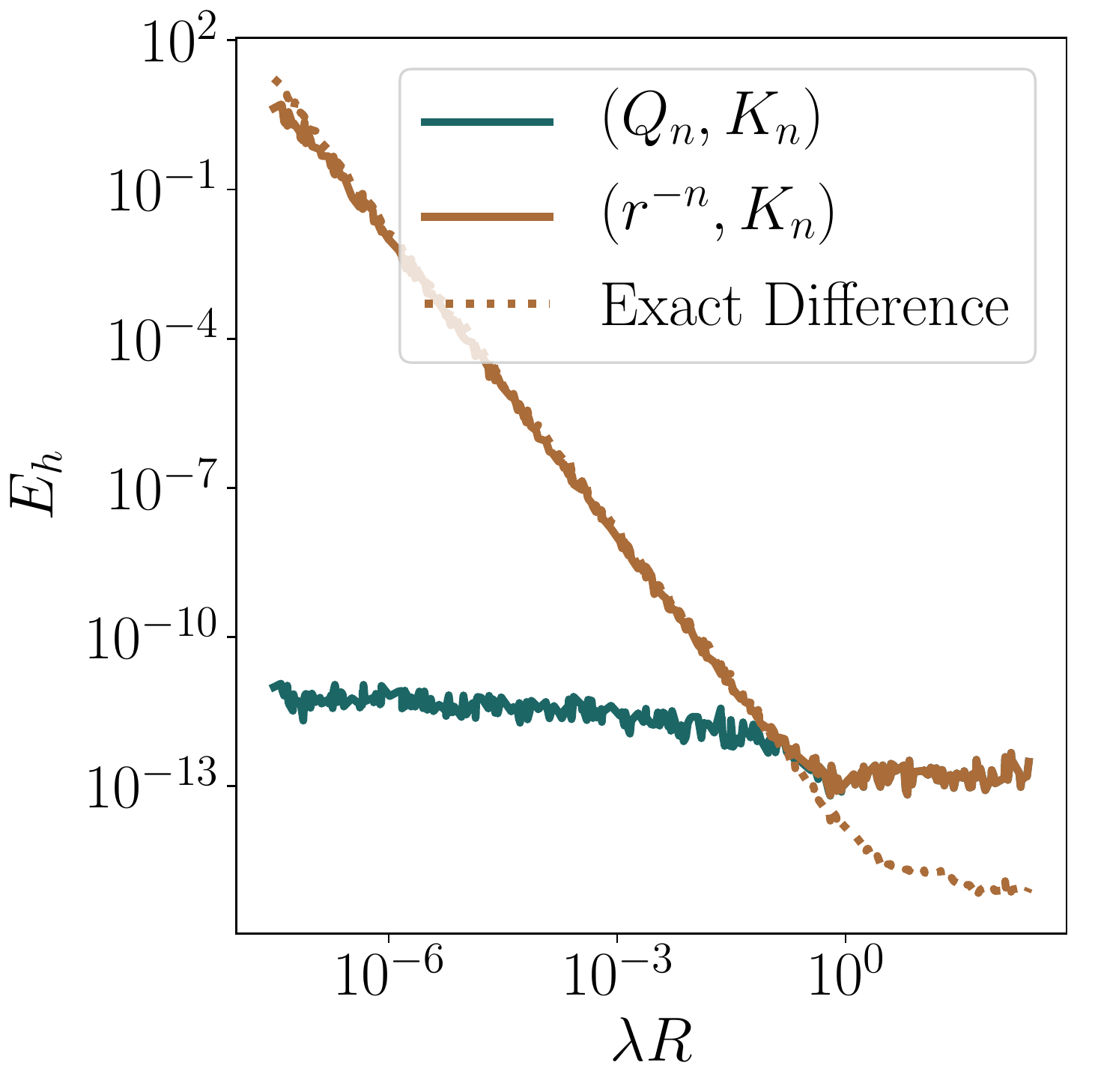}
    \caption{$E_h$ as $\lambda$ goes to zero.}
  \end{subfigure}

  \begin{subfigure}[t]{0.32\textwidth}
    \includegraphics[width=\linewidth]{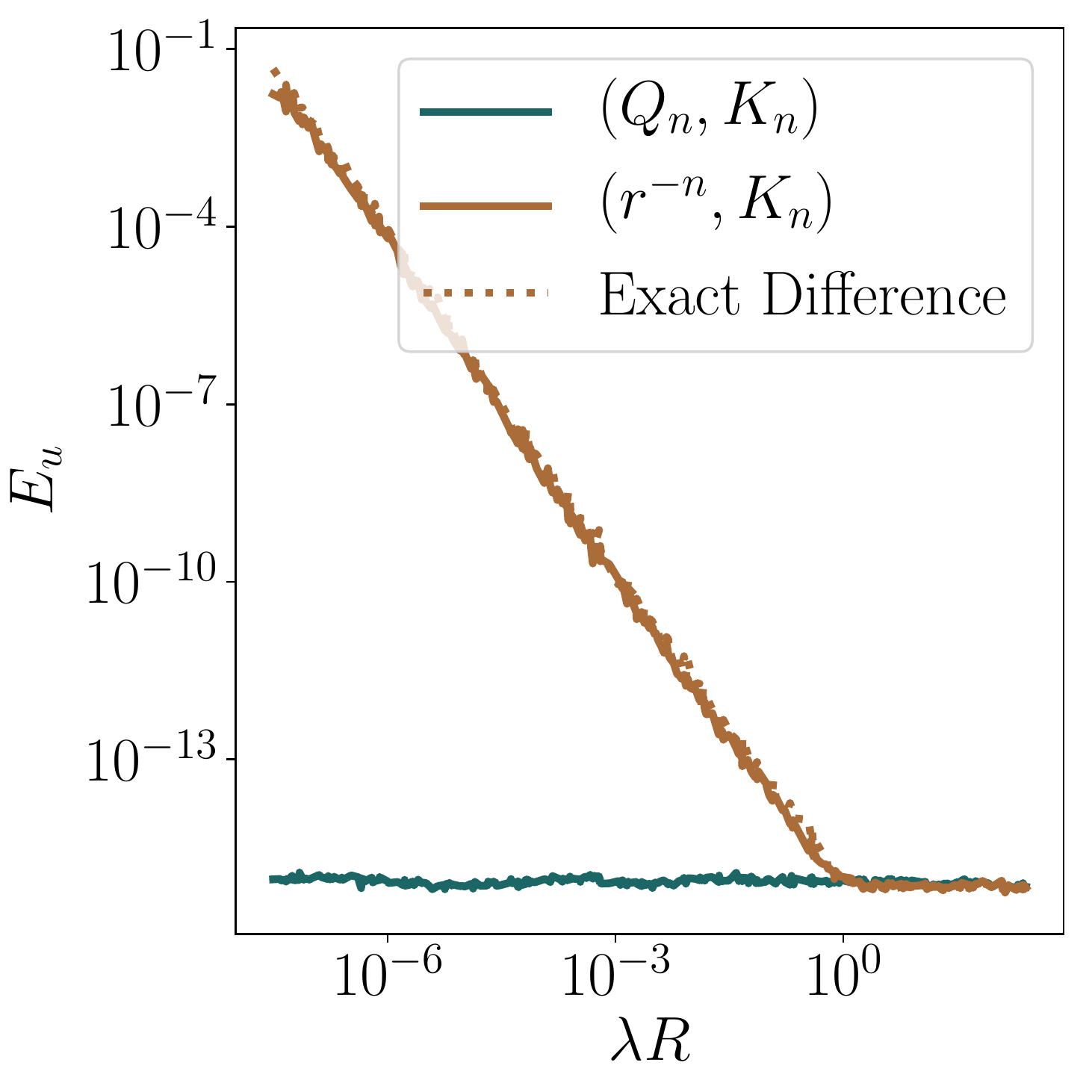}
    \caption{$E_u$ as $R$ goes to zero.}
  \end{subfigure}
  \begin{subfigure}[t]{0.32\textwidth}
    \includegraphics[width=\linewidth]{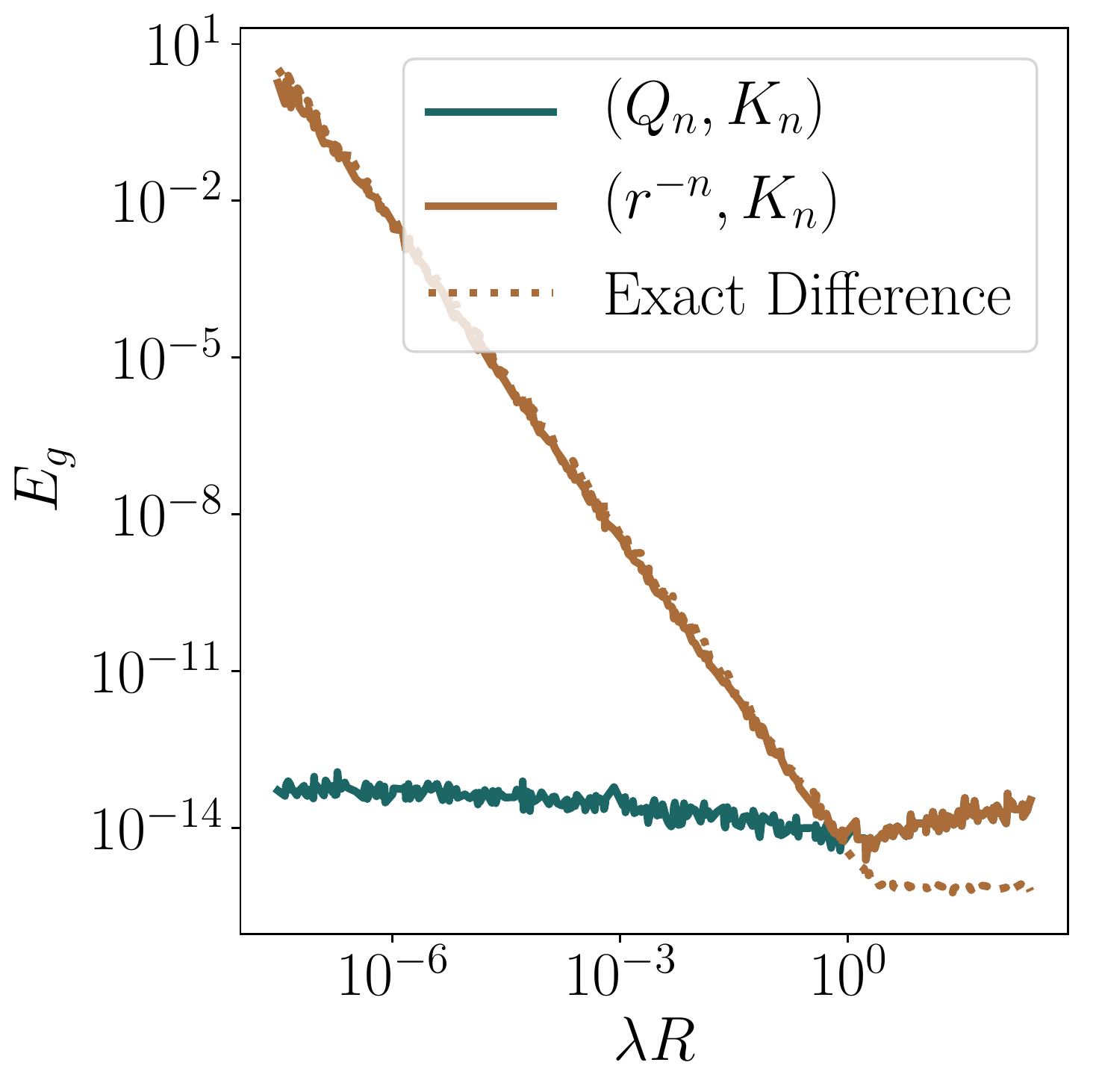}
    \caption{$E_g$ as $R$ goes to zero.}    
  \end{subfigure}
  \begin{subfigure}[t]{0.32\textwidth}
    \includegraphics[width=\linewidth]{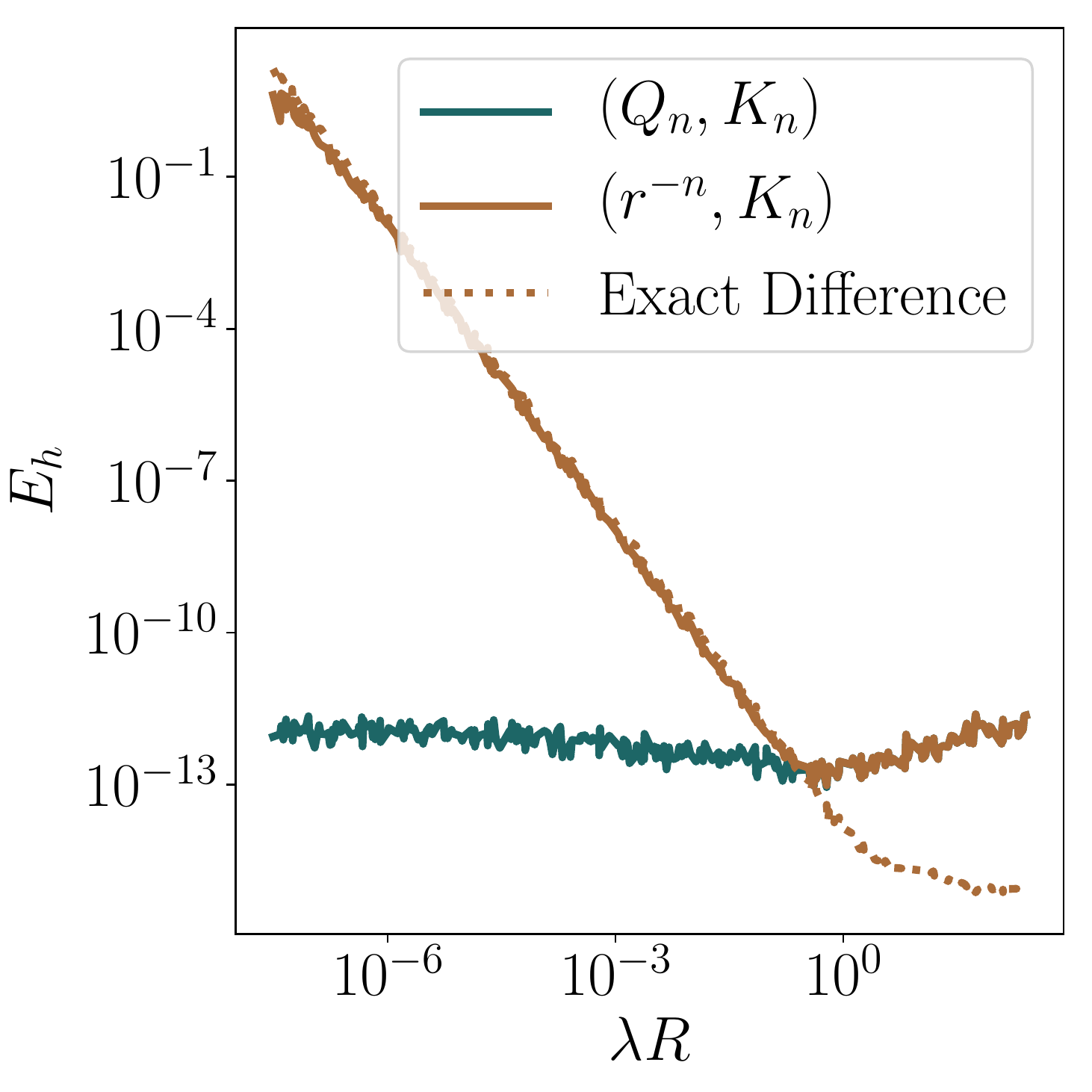}
    \caption{$E_h$ as $R$ goes to zero.}
  \end{subfigure}

  \caption{Exterior problem. In the top row, we plot
    the error measures
    as functions of $\lambda R$ for $R=0.5$ as $\lambda$
    goes to zero. In the bottom row, we plot the error
    measures as functions of $\lambda R$ for $\lambda = 0.5$
    as $R$ goes to zero.}

  \label{fig:errext}
  
\end{figure}

In \cref{fig:errint}, we plot the error measures as
functions of $\lambda R$ for both the limit as $R$
goes to zero with $\lambda$ fixed and vice-versa,
for the interior problem.
The behavior of the two limits is similar, in contrast
with the condition numbers of the previous section.
This is an indication that $\lambda R$ is
the relevant figure for applications. We see that the
new basis $(r^{|n|},P_n)$ is able to achieve high accuracy,
even as $\lambda R$ tends to zero. We note that $E_u$
is around machine precision, and there is some precision
loss for the errors of the derivatives, $E_g$ and $E_h$.
When using the na\"{i}ve basis, $(r^{|n|},I_n)$, on
the other hand, there is significant loss of accuracy as
$\lambda R$ goes to zero. 
The error in evaluating the exact difference between
the harmonic and modified Helmholtz parts agrees
well with the na\"{i}ve basis $(r^{|n|},I_n)$ in the small
$\lambda R$ limit. This error for the exact difference
shows that there is a fundamental problem in using the
basis $(r^{|n|},I_n)$ to represent such solutions.
For large $\lambda R$, the exact difference
is capable of near machine precision in even the derivatives,
as the procedure does not really involve numerical
differentiation. We see similar behavior for the
exterior problem in \cref{fig:errext}.

\section{Discussion}
\label{sec:discuss}

In the preceding, we have shown that the new bases
$(r^{|n|},P_n)$ and $(Q_n,K_n)$ offer significant advantages
over the na\"{i}ve approach when solving the modified
biharmonic equation on a disk. We now show how these
functions can be used to efficiently solve the modified
biharmonic equation \cref{eq:modbh} on more complex
geometries.

\subsection{Fast sums}

Let $\x_i$ be a set of $N$ points in space. Suppose that
we would like to evaluate the sum 

\begin{equation}
  \label{eq:modbhsum}
  u(\x_i) = \sum_{j\neq i}^{N} \lambda^2 c_j \mathcal{G}(\x_i,\x_j)
  + \lambda d_j \partial_{v_{j,1}} \mathcal{G}(\x_i,\x_j) +
  q_j \partial_{v_{j,2}v_{j,3}} \mathcal{G}(\x_i,\x_j) \, ,
\end{equation}
for each $\x_i$ efficiently. Direct evaluation would
be an $\bigo (N^2)$ calculation. When performing the analogous
sum with the Laplace Green's function, the fast multipole method
\cite{greengard1987fast,carrier1988fast} provides
a stable $\bigo (N)$ algorithm. We do not seek
to review the fast multipole method here, but we note
that such a method depends on a few key parts: a formula
for representing the sum due to a localized subset of the
points (a multipole expansion), a formula for representing
the sum due to points outside of a disk (a local expansion),
formulas for translating between these
representations (translation operators), and a hierarchical
organization of the points
in space. For details, see
\cite{greengard1987fast,carrier1988fast}.

Following the results of the previous section, we see that
an expansion in terms of the functions $(Q_n,K_n)$ provides
a stable representation for the sum due to points contained
inside a disk when evaluated at points sufficiently far from
that disk. Similarly, an expansion in terms of the functions
$(r^{|n|},P_n)$ can be used to stably represent the sum due to
points located sufficiently far outside of a disk. These are then
our multipole and local expansions, respectively. Starting
with the formulas for translating multipole and local expansions
for the Laplace and modified Helmholtz Green's functions,
which are included in \cref{sec:prelim}, it is straightforward
to derive translation operators for the $(Q_n,K_n)$ and
$(r^{|n|},P_n)$ expansions. We provide these in \cref{sec:props}.
Therefore, the new basis functions derived in this paper
result in a stable fast multipole method for the modified
biharmonic equation. Such a method has been implemented
by the author and results will be reported at a later date.

\subsection{Integral equations}

Because \cref{eq:modbh} is a homogeneous, fourth order
equation, it is well suited to solution using an integral
equation formulation. We do not attempt a review of the
literature here but point to \cite{chapko1997rothe,
  greengard1998integral,chapko2001combination,
  biros2004embedded,Jiang2013} for some representative 
examples. In an integral equation method, the solution 
is represented by a layer potential with unknown densities
defined on the boundary $\Gamma$,
i.e.

\begin{equation}
  \label{eq:layerpotential}
  u(\x) = \int_\Gamma K_1(\x,\y) \sigma_1(\y) +
  K_2(\x,\y) \sigma_2(\y) \, dS(\y) \; ,
\end{equation}
where the kernels $K_1$ and $K_2$ are typically defined
in terms of directional derivatives of the free-space
Green's function $\mathcal{G}$. For example, in
\cite{Jiang2013}, the kernels are
$K_1 = - \mathcal{G}_{\nu\nu} + \mathcal{G}_{\tau\tau}$
and $K_2 = -2 \mathcal{G}_{\nu\nu\nu} + 3(\Delta -\lambda^2)G_\nu
+ 2\lambda^2 \mathcal{G}_{\nu}$, where $\nu$ and $\tau$
represent the normal and tangential directions at
$\y$, respectively.

Continuing the example, the authors of \cite{Jiang2013}
then impose gradient boundary conditions on $u$, i.e.
they set $\partial_{\nu_x} u = f$ and $\partial_{\tau_x} u= g$
for some functions $f$ and $g$ defined on the boundary,
where $\nu_x$ and $\tau_x$ denote the normal and tangential
directions at a point $\x$ on the boundary. These boundary
conditions are of physical interest because they correspond
to ``no-slip'' boundary conditions for a stream function
representation of a fluid flow. Plugging
the form \cref{eq:layerpotential} into the boundary conditions,
we obtain the integral equation

\begin{equation}
  \label{eq:inteq}
  \begin{pmatrix}
    D_{11}(\x) & D_{12}(\x) \\
    0 & D_{22}(\x)
  \end{pmatrix}
  \begin{pmatrix}
    \sigma_1(\x) \\
    \sigma_2(\x)
  \end{pmatrix}
  + \int_\Gamma
  \begin{pmatrix}
    K_{11}(\x,\y) & K_{12}(\x,\y) \\
    K_{21}(\x,\y) & K_{22}(\x,\y)
  \end{pmatrix}
  \begin{pmatrix}
    \sigma_1(\y) \\
    \sigma_2(\y)
  \end{pmatrix}
  \, dS(\y)
  = \begin{pmatrix}
    f(\x) \\
    g(\x)
  \end{pmatrix} \; ,
\end{equation}
where $K_{11} = \partial_{\nu_x} K_1$, $K_{12} = \partial_{\nu_x} K_2$,
$K_{21} = \partial_{\tau_x} K_1$, and $K_{22} = \partial_{\tau_x} K_2$.
See \cite{Jiang2013} for details, including a simple preconditioner
for turning \cref{eq:inteq} into a well-conditioned second kind
integral equation and explicit formulas for the kernels $K_{ij}$.

In a Nystr\"{o}m discretization of the integral equation
\cref{eq:inteq}, the solution is represented by its
values at points on the curve $\Gamma$ which are used for
an integration rule. 
The basis of a Nystr\"{o}m method is then a
numerical quadrature of the integral
operator on the curve $\Gamma$. Typically, different
integration rules are required for
smooth, weakly singular, 
and singular integral kernels. Fortunately,
a number of quadrature rules
are available for handling these singularities
with high order accuracy
\cite{alpert,kapur,bremer2010,helsing08}.

Let $\x_i$
denote the points on $\Gamma$ of the Nystr\"{o}m
discretization, $\sigma_{1i} = \sigma_1(\x_i)$,
$\sigma_{2i} = \sigma_2(\x_i)$, $f_i = f(\x_i)$,
and $g_i = g(\x_i)$. In a slight abuse of
notation, the quadrature rule provides weights
$w_{ij}$ such that 

\begin{equation}
  \label{eq:inteqdisc}
  \begin{pmatrix}
    D_{11}(\x_i) & D_{12}(\x_i) \\
    0 & D_{22}(\x_i)
  \end{pmatrix}
  \begin{pmatrix}
    \sigma_{1i} \\
    \sigma_{2i}
  \end{pmatrix}
  + \sum_j w_{ij}
  \begin{pmatrix}
    K_{11}(\x_i,\x_j) & K_{12}(\x_i,\x_j) \\
    K_{21}(\x_i,\x_j) & K_{22}(\x_i,\x_j)
  \end{pmatrix}
  \begin{pmatrix}
    \sigma_{1j} \\
    \sigma_{2j}
  \end{pmatrix}
  = \begin{pmatrix}
    f_i \\
    g_i
  \end{pmatrix} \; ,
\end{equation}
is an accurate approximation of the original integral
equation at each point $\x_i$. The above is an abuse
of notation because the kernels are not often defined
when $\x = \y$, so that the formula is in general a
function of the kernel and the boundary. A key feature
of these integral rules is that the weight
$w_{ij}$ is typically a function of $j$ alone for
points $\x_i$ and $\x_j$ which are sufficiently far
apart. Therefore, much of the sum \cref{eq:inteqdisc}
is of the form \cref{eq:modbhsum} and the fast multipole
method described in the previous section can be used
to apply the operator on the left-hand-side
of \cref{eq:inteqdisc} rapidly.
Combined with an iterative solver such as GMRES
\cite{saad1986gmres}, this
provides a fast solution method for the densities
$\sigma_1$ and $\sigma_2$. The fast multipole
method can also be utilized to efficiently
evaluate the formula for $u$, \cref{eq:layerpotential},
at points inside the domain
\cite{helsing08,Kloeckner13,rachh2017fast}.

\section{Conclusion}

We have presented new special functions for representing
solutions of the modified biharmonic equation on both
the inside and outside of a disk. Numerical experiments
demonstrate the superiority of using these functions
over the na\"{i}ve approach, which would employ more familiar
functions associated with the Laplace and modified Helmholtz
equations. Further, we show how such functions can be used
to aid in the solution of the modified biharmonic equation
on more complex geometries with an integral equation approach.
We also present, in \cref{sec:props}, 
translation operators for multipole and local expansions
using our new radial basis functions. These are key
components of fast multipole-like methods 
\cite{carrier1988fast,greengard1987fast}.

The basic approach of the present paper applies to a number of
other Green's functions for high-order PDEs, 
assuming they can be expressed as the difference of 
Green's functions for lower order PDEs. This includes, for example,
the Green's function for the bi-Helmholtz equation 
\eqref{eq:bihelm}. Since the essential analysis concerns radial
functions only, the approach extends naturally to three dimensions.
We will explore these problems in subsequent work.



\section{Acknowledgments}

T.~Askham would like to thank Leslie Greengard
for many useful discussions. This work was supported
by the Office of the Assistant Secretary of Defense
for Research and Engineering and
AFOSR under NSSEFF Program Award FA9550-10-1-0180 and
award FA95550-15-1-0385.

\appendix
\gdef\thesection{\Alph{section}}
\makeatletter
\renewcommand\@seccntformat[1]{Appendix \csname the#1\endcsname.\hspace{0.5em}}
\makeatother

\section{Analytical preliminaries}
\label{sec:prelim}

For the Laplace kernel, the fast multipole
method is based on the manipulation of multipole
and Taylor expansions. The center of a multipole
expansion may be shifted using this formula:

\begin{lemma} \label{lemma:mpmp} 
(Adapted from Lemma 2.3 of \cite{greengard1987fast}.)
Suppose that

\begin{equation}
  \phi(z) = a_0 \log(z-z_0) + \sum_{l=1}^\infty
\dfrac{a_l}{(z-z_0)^l}  
\end{equation}
is a multipole expansion of the potential due
to a charge density which is contained inside
the disk of radius $R$ about $z_0$. Then, for 
$z$ outside the disk of radius $R+|z_0|$ about
the origin

\begin{equation}
  \phi(z) = b_0 \log(z) + \sum_{l=1}^\infty
\dfrac{b_l}{z^l}  \; ,
\end{equation}
where $b_0 = a_0$ and 

\begin{equation}
  b_l = \left ( \sum_{m=1}^l a_m z_0^{l-m} {l-1 \choose m-1} \right )
  - \dfrac{a_0 z_0^l}{l} \; ,
\end{equation}
using the standard notation for binomial 
coefficients. We also have the following bound
for the truncation error. With $p \geq 1$,

\begin{equation}
  \left | \phi(z) - b_0 \log(z) - 
\sum_{l=1}^p \dfrac{b_l}{z^l}  \right |
\leq \dfrac{F/2\pi}{1-(|z_0| + R)/|z|}
\left ( \dfrac{|z_0|+R}{|z|} \right )^{p+1} \; ,
\end{equation}
where $F$ is the $L_1$ norm of the density.
\end{lemma}

A multipole expansion may be converted into a Taylor
expansion using the following formula:
\begin{lemma} \label{lemma:mploc} 
(Adapted from Lemma 2.4 of \cite{greengard1987fast}.)
Suppose that a charge density is contained inside
the disk of radius $R$ about $z_0$ with $|z_0| > (1+c)R$
for some $c>1$. Let the multipole expansion due
to this density be given as in Lemma \ref{lemma:mpmp}. 
Then, this multipole expansion converges inside 
the disk of radius $R$ about the origin and 
can be represented by a power series there:

\begin{equation}
  \phi(z) = \sum_{l=1}^\infty b_l z^l
\end{equation}
where,

\begin{equation}
  b_0 = \sum_{m=1}^\infty \dfrac{a_m}{z_0^m} \left(-1 \right)^m
+ a_0 \log(-z0) \; ,
\end{equation}
and

\begin{equation}
  b_l = \dfrac{1}{z_0^l} \left (\sum_{m=1}^\infty \dfrac{a_m}{z_0^m}
{l+m-1 \choose m-1} \left (-1 \right)^m \right )
 - \dfrac{a_0}{z_0^l \ l} \; ,
\end{equation}
for $l \geq 1$. There is a similar error bound for this lemma, see
\cite{greengard1987fast} for details. Suppose that  the charge density is supported 
in a box and the evaluation points $z$ are taken in another. 
In the case that
these two boxes are well separated from each other, we have

\begin{equation}
\left | \phi(z) - \sum_{l=0}^p b_l z^l \right | \leq CF 
\left ( \dfrac12 \right )^p \; ,
\end{equation}
where $F$ is as in Lemma \ref{lemma:mpmp}
and $C$ is a constant.
\end{lemma}

The center of a Taylor expansion can be shifted
using the following formula:

\begin{lemma}
\label{lemma:locloc} (Adapted from Lemma 2.5 in \cite{greengard1987fast})
Let $z_0$, $z$, and $a_l$ for $l = 0,1,\ldots, p$ be complex. 
Then

\begin{equation}
  \sum_{l=0}^p a_l (z-z_0)^l = \sum_{m=0}^p 
\left ( \sum_{l=m}^p a_l {l \choose m} (-z_0)^{l-m} \right)
z^m \; .
\end{equation}
This formula is exact.
\end{lemma}

The fast multipole method for the modified Helmholtz
kernel is based on the manipulation of expansions
in the Bessel functions $I_n$ and $K_n$. The formula
for shifting the center of an expansion in the $K_n$
functions is:

\begin{lemma} (Adapted from \cite{cheng2006adaptive}.) Suppose that

\begin{equation}
\phi(\x) = \sum_{l = -\infty}^\infty a_l  
K_l(\lambda \rho') e^{il\theta'} \;,
\end{equation}
where $(\rho',\theta')$ are the polar coordinates of $\x$
with respect to the point $\x_0$, 
is a multipole expansion of the potential due to 
a charge density which is contained inside the 
disk of radius $R$ about $\x_0$. Let $(\rho_0,\theta_0)$
be the polar coordinates of $\x_0$ with respect to the
origin. Then, for $\x$ outside the disk of radius 
$R + \rho_0$, we have

\begin{equation}
  \phi(\x) = \sum_{l=-\infty}^\infty b_l K_l(\lambda \rho) e^{il\theta}
  \; ,
\end{equation}
where $(\rho,\theta)$ are the coordinates of $\x$ with
respect to the origin and the translated coefficients
are given by

\begin{equation}
  b_l = \sum_{m=-\infty}^\infty a_m I_{l-m}(\lambda \rho_0)
  e^{-i(l-m)\theta_0} \; .
\end{equation}
  
\end{lemma}

An expansion in the $K_n$ functions may be converted
into an expansion in the $I_n$ functions using the
formula:

\begin{lemma} (Adapted from \cite{cheng2006adaptive}.) Suppose that

\begin{equation}
\phi(\x) = \sum_{l = -\infty}^\infty a_l  
K_l(\lambda \rho') e^{il\theta'} \;,
\end{equation}
where $(\rho',\theta')$ are the polar coordinates of $\x$
with respect to the point $\x_0$, 
is a multipole expansion of the potential due to 
a charge density which is contained inside the 
disk of radius $R < |\x_0|$ about $\x_0$. Let $(\rho_0,\theta_0)$
be the polar coordinates of $\x_0$ with respect to the
origin. Then, for $\x$ within the disk of radius 
$\rho_0 - R$, we have

\begin{equation}
  \phi(\x) = \sum_{l=-\infty}^\infty b_l I_l(\lambda \rho) e^{il\theta}
  \; ,
\end{equation}
where $(\rho,\theta)$ are the coordinates of $\x$ with
respect to the origin and the translated coefficients
are given by

\begin{equation}
  b_l = \sum_{m=-\infty}^\infty a_m (-1)^m K_{l-m}(\lambda \rho_0)
  e^{-i(l-m)\theta_0} \; .
\end{equation}
\end{lemma}

We phrase the following translation operator in terms of
a parent and child box. A child box is any of the four
boxes resulting from dividing a parent box into equal
square quadrants.
The center of an expansion in the $I_n$ functions
may be shifted using the following formula:

\begin{lemma} (Adapted from \cite{cheng2006adaptive}.) Suppose that

\begin{equation}
\Psi(\x) = \sum_{l = -\infty}^\infty a_l  
I_l(\lambda \rho) e^{il\theta} \;,
\end{equation}
where $(\rho,\theta)$ are the polar coordinates of $\x$
with respect to the origin, 
is a local expansion for a parent box centered at
the origin. Let $(\rho_0,\theta_0)$
be the polar coordinates of the center $\x_0$ 
of a child box with respect to the origin. 
Then, for $\x$ inside the child box, we have

\begin{equation}
  \Psi(\x) = \sum_{l=-\infty}^\infty b_l I_l(\lambda \rho') e^{il\theta'}
  \; ,
\end{equation}
where $(\rho',\theta')$ are the coordinates of $\x$ with
respect to the child's center $\x_0$ and the new local 
expansion coefficients are given by

\begin{equation}
  b_l = \sum_{m=-\infty}^\infty a_m I_{l-m}(\lambda \rho_0)
  e^{-i(l-m)\theta_0} \; .
\end{equation}
  
\end{lemma}

\section{Translation properties of $Q_n$ and $P_n$}
\label{sec:props}

The stabilized fast multipole method for the modified
biharmonic equation that we develop is based on
manipulating multipole expansions in the functions
$(Q_n,K_n)$ and local expansions in the functions
$(r^{|n|},P_n)$. The center of such a multipole expansion
can be shifted using the following formula:

\begin{lemma}
  Suppose that 
\begin{equation}
\phi(\x) = \sum_{l = -\infty}^\infty \left (a_l Q_l(\rho')
+ b_l K_l(\lambda \rho') \right ) e^{il\theta'} \;,
\end{equation}
where $(\rho',\theta')$ are the polar coordinates of $\x$
with respect to the point $\x_0$, 
is a multipole expansion of the potential due to 
a charge density which is contained inside the 
disk of radius $R$ about $\x_0$. Let $(\rho_0,\theta_0)$
be the polar coordinates of $\x_0$ with respect to the
origin. Then, for $\x$ outside the disk of radius 
$R + \rho_0$, we have

\begin{equation}
  \phi(\x) = \sum_{l = -\infty}^\infty \left (c_l Q_l(\rho)
+ d_l K_l(\lambda \rho) \right ) e^{il\theta} \;,
\end{equation}
where $(\rho,\theta)$ are the coordinates of $\x$ with
respect to the origin. We have that $c_0 = a_0$ and 

\begin{align}
d_0 &= \sum_{m=-\infty}^{-1} (a_m+b_m) I_{-m}(\lambda \rho_0)  
 e^{-i(l-m)\theta_0} + a_0 P_0(\lambda \rho_0) 
 + b_0 I_0(\lambda \rho_0) \nonumber \\
 & \quad + \sum_{m=1}^{\infty} (a_m+b_m) I_{-m}(\lambda \rho_0)
  e^{-i(l-m)\theta_0} \; .
\end{align}
For $l > 0$ the translated coefficients
are given by

\begin{equation}
  c_l = \sum_{m=0}^l a_m \left (\dfrac{\lambda \rho_0}{2} \right)^{l-m}
  \dfrac{e^{-i(l-m)\theta_0}}{(l-m)!}
\end{equation}
and
\begin{align}
  d_l &= \sum_{m=-\infty}^{-1} (a_m+b_m) I_{l-m}(\lambda \rho_0)
  e^{-i(l-m)\theta_0} \nonumber \\ 
  & \qquad + \sum_{m=0}^l (a_m P_{l-m}(\rho_0) + b_m I_{l-m}(\lambda \rho_0))
  e^{-i(l-m)\theta_0} \nonumber \\
  & \qquad + \sum_{m=l+1}^{\infty} (a_m+b_m) I_{l-m}(\lambda \rho_0)
  e^{-i(l-m)\theta_0} \; .
\end{align}
For $l < 0$ the translated coefficients
are given by

\begin{equation}
  c_l = \sum_{m=l}^0 a_m \left (\dfrac{\lambda \rho_0}{2} \right)^{|l-m|}
  \dfrac{e^{-i(l-m)\theta_0}}{|l-m|!}
\end{equation}
and
\begin{align}
  d_l &= \sum_{m=-\infty}^{l-1} (a_m+b_m) I_{l-m}(\lambda \rho_0)
  e^{-i(l-m)\theta_0} \nonumber \\ 
  & \qquad + \sum_{m=l}^0 (a_m P_{l-m}(\rho_0) + b_m I_{l-m}(\lambda \rho_0))
  e^{-i(l-m)\theta_0} \nonumber \\
  & \qquad + \sum_{m=1}^{\infty} (a_m+b_m) I_{l-m}(\lambda \rho_0)
  e^{-i(l-m)\theta_0} \; .
\end{align}
  
\end{lemma}

A multipole expansion can be converted into a local expansion
using the following formula:

\begin{lemma}
  Suppose that 
\begin{equation}
\phi(\x) = \sum_{l = -\infty}^\infty \left (a_l Q_l(\rho')
+ b_l K_l(\lambda \rho') \right ) e^{il\theta'} \;,
\end{equation}
where $(\rho',\theta')$ are the polar coordinates of $\x$
with respect to the point $\x_0$, 
is a multipole expansion of the potential due to 
a charge density which is contained inside the 
disk of radius $R < |\x_0|$ about $\x_0$. Let $(\rho_0,\theta_0)$
be the polar coordinates of $\x_0$ with respect to the
origin. Then, for $\x$ within the disk of radius 
$\rho_0-R$, we have

\begin{equation}
  \phi(\x) = \sum_{l=-\infty}^{\infty} \left ( c_l P_l(\lambda \rho) +
  d_l (\lambda \rho)^{|l|}  \right ) e^{il\theta} \; ,
\end{equation}
where $(\rho,\theta)$ are the coordinates of $\x$ with
respect to the origin. For all $l$, we have that

\begin{equation}
  c_l = \sum_{m=-\infty}^\infty (-1)^m (a_m+b_m) 
  K_{l-m}(\lambda \rho_0) e^{i(l-m)\theta_0} \; .
\end{equation}
For $l > 0$, we have

\begin{align}
  d_l &= \dfrac{1}{2^ll!} 
  \sum_{m=-\infty}^{0} (-1)^m (a_m Q_{l-m}(\rho_0) 
  + b_m K_{l-m}(\lambda \rho_0)) e^{-i(l-m)\theta_0}  \nonumber \\
  & \quad + \dfrac{1}{2^ll!} \sum_{m=1}^\infty (-1)^m
  (a_m+b_m) K_{l-m}(\lambda \rho_0) e^{-i(l-m)\theta_0} \; .
\end{align}
For $l = 0$, we have

\begin{align}
  d_0 &= \sum_{m=-\infty}^{\infty} (-1)^m (a_m Q_{-m}(\rho_0) 
  + b_m K_{-m}(\lambda \rho_0)) e^{im\theta_0} 
\end{align}
For $l < 0$, we have

\begin{align}
  d_l &= \dfrac{1}{2^{|l|}|l|!} 
  \sum_{m=-\infty}^{-1} (-1)^m (a_m + b_m) K_{l-m}(\lambda \rho_0) 
  e^{-i(l-m)\theta_0}  \nonumber \\
  & \quad + \dfrac{1}{2^{|l|}|l|!} \sum_{m=0}^\infty (-1)^m
  (a_mQ_{l-m}(\rho_0)+b_mK_{l-m}(\lambda \rho_0))e^{-i(l-m)\theta_0} \; .
\end{align}

\end{lemma}

Finally, we can shift the center of a local expansion in
the $(r^{|n|},P_n)$ using the following formula:

\begin{lemma} Suppose that

\begin{equation}
\Psi(\x) = \sum_{l = -\infty}^\infty (a_l  
P_l(\rho) + b_l (\lambda \rho)^{|l|} ) e^{il\theta} \;,
\end{equation}
where $(\rho,\theta)$ are the polar coordinates of $\x$
with respect to the origin, 
is a local expansion for a parent box centered at
the origin. Let $(\rho_0,\theta_0)$
be the polar coordinates of the center $\x_0$ 
of a child box with respect to the origin. 
Then, for $\x$ inside the child box, we have

\begin{equation}
\Psi(\x) = \sum_{l = -\infty}^\infty (c_l  
P_l(\rho') + d_l (\lambda \rho')^{|l|} ) e^{il\theta'} \;,
\end{equation}
where $(\rho',\theta')$ are the coordinates of $\x$ with
respect to the child's center $\x_0$. For all $l$, we have

\begin{equation}
  c_l = \sum_{m=-\infty}^\infty a_m I_{l-m}(\lambda \rho_0)
  e^{-i(l-m)\theta_0} \; .
\end{equation}
For $l > 0$ , we have 

\begin{align}
  d_l &= \dfrac{1}{2^ll!} \sum_{m=-\infty}^{l-1} a_m I_{l-m}(\lambda \rho_0)
  e^{-i(l-m)\theta_0} \nonumber \\
  & \quad + \dfrac{1}{2^ll!} \sum_{m=l}^{\infty} a_m P_{l-m}(\lambda \rho_0)
  e^{-i(l-m)\theta_0} \nonumber \\
  & \quad + \sum_{m=l}^\infty {m\choose l} b_m \rho_0^{|l-m|} e^{-i(l-m)\theta_0}
  \; .
\end{align}
For $l = 0$ , we have 

\begin{align}
  d_0 &= \sum_{m=-\infty}^{\infty} a_m P_{-m}(\lambda \rho_0)
  e^{im\theta_0} + \sum_{m=-\infty}^\infty b_m \rho_0^{|m|} e^{im\theta_0}
  \; .
\end{align}
For $l < 0$ , we have 

\begin{align}
  d_l &= \dfrac{1}{2^{|l|}|l|!} \sum_{m=-\infty}^{l} a_m P_{l-m}(\lambda \rho_0)
  e^{-i(l-m)\theta_0} \nonumber \\
  & \quad + \dfrac{1}{2^{|l|}|l|!} 
  \sum_{m=l+1}^{\infty} a_m I_{l-m}(\lambda \rho_0)
  e^{-i(l-m)\theta_0} \nonumber \\
  & \quad + \sum_{m=-\infty}^l {|m|\choose |l|} b_m 
  \rho_0^{|l-m|} e^{-i(l-m)\theta_0} \; .
\end{align}
  
\end{lemma}

\bibliographystyle{elsarticle-num}
\bibliography{refs}

\end{document}